\documentstyle[psfig]{article}
\newtheorem{theorem}{Theorem}
\newtheorem{assumption}{Assumption}

\def\endproof{\vrule height7pt width7pt depth0pt}

\newcommand{\carre}[2]{\put(#1){\framebox(10,10){$#2$}}}

\newcommand{\rond}[2]{\put(#1){\makebox(0,0){$#2$}}%
    \put(#1){\circle{12}}}

\newcommand{\coinbl}[2]{\put(#1){\makebox(0,0)[bl]{$#2$}}}
\newcommand{\coinbr}[2]{\put(#1){\makebox(0,0)[br]{$#2$}}}
\newcommand{\cointl}[2]{\put(#1){\makebox(0,0)[tl]{$#2$}}}
\newcommand{\cointr}[2]{\put(#1){\makebox(0,0)[tr]{$#2$}}}

\title{A mechanical model \\
for the transportation problem}
\author{M. H\'enon \\
C.N.R.S., Observatoire de Nice, \\
BP 4229, 06304 Nice Cedex 4, France}
\date{May 18, 1992}
		
\begin{document}
\maketitle

\begin{abstract}

We describe a mechanical device which can be used as an analog
computer to solve the transportation problem. In practice
this device is simulated by a numerical algorithm. Tests show that
this algorithm is 60 times faster than a current subroutine
(NAG library) for an average $1000 \times 1000$ problem. Its
performance is even better for degenerate problems in which the
weights take only a small number of integer values.

\vspace{12pt}

\end{abstract}

\newpage
\tableofcontents
\newpage

\section{Introduction}

We describe here an algorithm for the solution of the
transportation problem \cite{Hit41a} (also known as the Hitchcock
problem).

The development of this algorithm had its origin in studies of the
lattice gas method for three-dimensional fluid simulations
\cite{HLF86a}.  The optimization of the collision table has
generally the form of a transportation problem \cite{Hen89a,RHF*88a},
with large cost matrices.  Classical algorithms were found to
require prohibitively long computing times. Therefore an attempt
was made to devise a method which would take advantage of the
peculiarities of the lattice gas problem. This method then turned
out to be of general applicability.

The present algorithm was developed independently of the already
published studies of the transportation problem and related
optimization problems.  This was not planned; it only reflects the
way things happened, and the ignorance of this author who comes
from a rather different field.  More will be said about this in
Section~\ref{s:conclusions}.

The paper is organized as follows. Section~\ref{s:problem}
defines the problem. In Section~\ref{s:analog}, we
describe a mechanical device which can be used as an analog
computer to solve the transportation problem.  In
Section~\ref{s:num-sim} we develop an appropriate graph
representation and a numerical scheme which simulates the
mechanical model. This is illustrated by a detailed example in
Section~\ref{s:example}. In Section~\ref{s:formal}, we give a
rigorous definition and justification of the algorithm.
Section~\ref{s:implementation} describes some aspects of the
computer implementation. In Section~\ref{s:tests}, the algorithm
is compared with the NAG library subroutine for the solution of
the transportation problem. In the particular case of the assignment
problem, comparisons are also made with the algorithm of Burkard
and Derigs \cite{BD80a}. A few comments are made in
Section~\ref{s:conclusions}. Finally, an Appendix derives some
bounds on the number of operations.
\section{The problem}
\label{s:problem}

We are given a matrix $c_{ij}$ and two vectors $a_i \ge 0$, $b_j
\ge 0$. The index $i$ runs from 1 to $m$ and the index $j$ runs
from 1 to $n$.  There is

		\begin{equation}
\sum _{i=1}^m a_i = \sum_{j=1}^n b_j.
					\label{zerosum}
		\end{equation}
The problem is to find coefficients $f_{ij}$ which {\em maximize}
the sum
		\begin{equation}
\sum_i \sum_j c_{ij} f_{ij}
					\label{maxisum}
		\end{equation}
subject to the constraints
		\begin{eqnarray}
\sum_j f_{ij} & = & a_i \qquad (i = 1, \ldots, m), 
					\label{rowsum} \\
\sum_i f_{ij} & = & b_j \qquad (j = 1, \ldots, n), 
					\label{columnsum} \\
f_{ij} & \ge & 0.
					\label{positivef}
		\end{eqnarray}
A set of $f_{ij}$ which satisfies the constraints (\ref{rowsum})
to (\ref{positivef}) will be called a {\em feasible solution}. A
feasible solution which maximizes (\ref{maxisum}) will be called
an {\em optimal solution}.

Notes: (i) This is essentially the {\em transportation problem}.  It
can be reduced to the standard form given for instance in
\cite[Section~7.4]{PS82a}, by defining
		\begin{equation}
c^*_{ij} = c_{\rm sup} - c_{ij},
		\end{equation}
where $c_{\rm sup}$ is a constant satisfying
		\begin{equation}
c_{\rm sup} > \max_{i,j} c_{ij}.
					\label{defcsup}
		\end{equation}

(ii) There is no sign condition on the $c_{ij}$, which can be 
positive or negative.

(iii) If $a_i = 0$ for one row $i$, then from (\ref{rowsum}) and
(\ref{positivef}) we have $f_{ij} = 0$ for all $j$.  This row does
not contribute to the sum (\ref{maxisum}); the values of the
$c_{ij}$ on that row are irrelevant.  Thus this row could be
eliminated without changing the problem.  The same holds if $b_j =
0$ for some $j$.  We could therefore in principle restrict our
attention to the case where the $a_i$ and $b_j$ are strictly
positive, as is usually done \cite{Had62a,PS82a}. In practice,
however, it is convenient to be able to include the cases with
some $a_i = 0$ and/or $b_j = 0$ into the general treatment. The
algorithm to be described works just as well in such cases.

(iv) The $a_i$, $b_j$, $c_{ij}$ can be integer or real numbers.

(v) A special case of the transportation problem is
		\begin{equation}
a_i = 1 \qquad (i = 1, \ldots, m), 
\qquad b_j = 1 \qquad (j = 1, \ldots, n).
		\end{equation}
 From (\ref{zerosum}) it follows then that
		\begin{equation}
m = n.
		\end{equation}

If we prescribe the additional constraint 
		\begin{equation}
f_{ij} = 0 \mbox{ or } 1,
					\label{0or1}
		\end{equation}
we obtain the classical {\em assignment problem}
(\cite[chap.~11]{PS82a}).
\section{A mechanical model}
\label{s:analog}

We show now that it is possible to build a simple mechanical
device which solves the transportation problem. This device acts as a
{\em analog computer}: the numbers entering the problem are
represented by physical quantities, and the equations are
replaced by physical laws. 

\begin{figure}[hbtp]
\centerline{\psfig{file=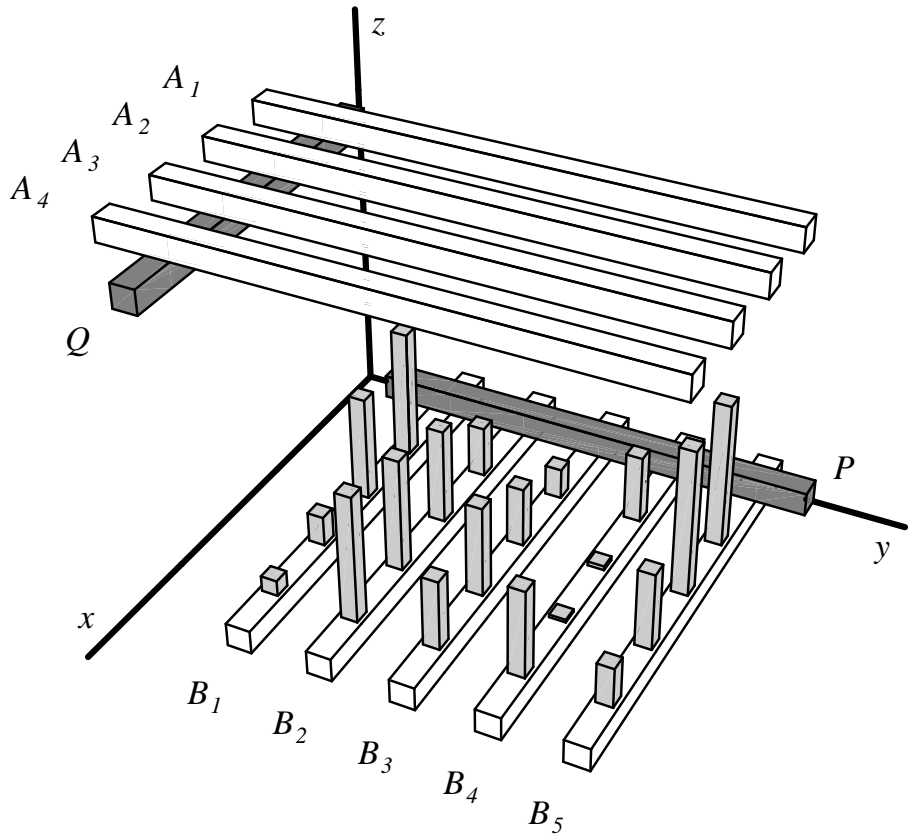,width=12cm,%
      bbllx=32mm,bblly=73mm,bburx=170mm,bbury=199mm}}
\caption{
An analog computer for the solution of the transportation problem.}
\label{f:modele}
		\end{figure}

We define a system of axes $x$, $y$, $z$ in physical space
(Fig.~\ref{f:modele}).

To every value of $i$ is associated a rod $A_i$ parallel to the
$y$ axis, which we will call a {\em row} by reference to the
$c_{ij}$ matrix. Mechanical constraints (not shown on the figure)
ensure that each row can only move in the vertical direction. More
precisely, each row remains parallel to the $y$ axis and moves in
a fixed vertical plane $x = $ Const. The variable height $z$ of
the lower face of row $A_i$ will be designated by $\alpha_i$.  Row
$A_i$ has a weight $a_i$ and thus is subjected to a force $a_i$
towards the negative $z$ axis.

Similarly, to every value of $j$ is associated a rod $B_j$
parallel to the $x$ axis, which we will call a {\em column}.
(Notice that these ``columns'' are horizontal in the
three-dimensional physical space !). Each column is constrained to
move only vertically: it remains parallel to the $x$ axis and
moves in a fixed vertical plane $y = $ Const. It has a {\em
negative} weight $- b_j$ (or, if one prefers, a buoyancy $b_j$)
and thus is subjected to a force $b_j$ towards the positive $z$
axis. (In an actual model, this might be realized with cables and
counterweights).  We call $z = \beta_j$ the height of its upper
face.

Finally, small vertical cylinders or {\sl studs} of height
$c_{ij}$ and of negligible weight are placed on the columns, in
such a way that each stud enforces a minimal distance between row
$A_i$ and column $B_j$:

		\begin{equation}
\alpha_i - \beta_j \ge c_{ij}.
					\label{plot}
\end{equation} 

(Note: this description seems to imply that $c_{ij} \ge 0$.
Actually it is possible, although mechanically more awkward, to
have negative values of the $c_{ij}$ by bending the rods. One can
also make all $c_{ij}$ positive by adding a sufficiently large
constant to all of them.  Therefore we continue to consider that
the $c_{ij}$ can be arbitrary.)

The potential energy of the system is, within an additive
constant:
		\begin{equation}
U = \sum_i a_i \alpha_i - \sum_j b_j \beta_j.
		\end{equation}

Initially, all rods are maintained at a fixed position by two
additional fixed rods $P$ and $Q$ acting as {\sl stops}
(Fig.~\ref{f:modele}), with the rows $A_i$ well above the columns
$B_j$, so that there is no contact between the rows and the
studs.  For instance we take $\beta_j = 0$ ($j = 1, \ldots, n$)
and  $\alpha_i = c_{\rm sup}$ ($i = 1, \ldots, m$).  Then the rods
are released by removing the stops $P$ and $Q$, and the system
starts evolving. Rows go down, columns go up,
and contacts are made with the studs. Aggregates of rows and
columns are progressively formed. As new contacts are made, these
aggregates are modified. Thus a complex evolution may take
place.

 It will be convenient to imagine that the system is immersed in a
viscous fluid, so that the velocity of an object, rather than its
acceleration, is proportional to the force to which it is
subjected. More specifically, let us define formally an aggregate
as a subset $S$ of rows and columns which are in contact (they
form a connected set in space) and move with the same velocity. We
also assume that it is a maximal connected set: no other row or
column touches it. The downward force acting on this aggregate is
its total weight, which we call $M(S)$:

		\begin{equation}
M(S) = \sum_{A_i \in S} a_i - \sum_{B_j \in S} b_j.
					\label{ag-weight}
		\end{equation}
We will simply assume that the aggregate moves with a velocity
		\begin{equation}
{dz \over dt} = - M(S).
					\label{ag-motion}
		\end{equation}
In particular an aggregate remains motionless if the force
applied to it vanishes.

It is intuitively clear that the system will reach an equilibrium,
an in fact we have
		\begin{theorem}
An equilibrium is reached after a finite time.
		\end{theorem}
Proof: an aggregate $S$ has a potential energy
		\begin{equation}
U(S) = \sum_{A_i \in S} a_i \alpha_i 
  - \sum_{B_j \in S} b_j \beta_j.
		\end{equation}
 From (\ref{ag-weight}) and (\ref{ag-motion}) we find that this
potential energy decreases with time according to
		\begin{equation}
{d U(S) \over dt} = - M^2(S).
		\end{equation}
We call $|M|_{\rm min}$ the minimum of all non-zero values of
$|M(S)|$, over all subsets $S$ of the full set of rods. Since
there is only a finite number of subsets, we have $|M|_{\rm min} >
0$, and:
		\begin{equation}
\cases{& either $d U(S) / dt = 0$, \cr
       & or $d U(S) / dt \le - |M|^2_{\rm min}$. \cr}
		\end{equation}
The total potential energy $U$ is the sum of the potential
energies of the aggregates. Therefore we also have
		\begin{equation}
\cases{& either $d U / dt = 0$, \cr
       & or $d U / dt \le - |M|^2_{\rm min}$. \cr}
		\end{equation}
The first case is realized only if $d U(S) / dt = 0$ for every
aggregate, i.e. if all aggregates are motionless. Thus: either
the system is in equilibrium, or its potential energy decreases at
a rate at least equal to $|M|^2_{\rm min}$.

On the other hand we have
		\begin{equation}
\alpha_i - \beta_j \ge \min_{i,j} c_{ij}.
		\end{equation}
Multiplying by $a_i b_j$, summing on $i$ and $j$, and using
(\ref{zerosum}), we obtain
		\begin{equation}
U \sum a_i \ge \min_{i,j} c_{ij} \left( \sum a_i \right)^2.
		\end{equation}
This gives a lower bound for $U$. 

Combining these results, we find that the system reaches an
equilibrium after a finite time (for which an upper bound is
easily derived).~\endproof

We consider now such an equilibrium state. We will show that
		\begin{theorem}
\label{th:equiv}
If the system is in equilibrium, and if $F_{ij}$ is the force
transmitted through stud $c_{ij}$ from row $A_i$ to column $B_j$,
then $f_{ij} = F_{ij}$ is an optimal solution of the transportation
problem.
		\end{theorem}
Proof: (i) Each row is in equilibrium, therefore
		\begin{equation}
\sum_j F_{ij} = a_i.
					\label{rowequil}
		\end{equation}

(ii) Each column is in equilibrium, therefore
		\begin{equation}
\sum_i F_{ij} = b_j.
					\label{columnequil}
		\end{equation}

(iii) $F_{ij}$ cannot be negative since it is transmitted by
contact:
		\begin{equation}
F_{ij} \ge 0.
		\end{equation}
Therefore $F_{ij}$ is a feasible solution.

(iv) If $F_{ij} >0$, row $A_i$ is in contact with column $B_j$,
and therefore
		\begin{equation}
\alpha_i - \beta_j = c_{ij}.
		\end{equation}
It follows that
		\begin{equation}
F_{ij} (\alpha_i - \beta_j - c_{ij}) = 0 \qquad \forall i, j.
					\label{slackness}
		\end{equation}
Summing (\ref{slackness}) over $i$ and $j$ and using
(\ref{rowequil}) and (\ref{columnequil}), we obtain
		\begin{equation}
\sum_i a_i \alpha_i - \sum_j b_j \beta_j 
- \sum_i \sum_j F_{ij} c_{ij} = 0.
					\label{bidon1}
		\end{equation}
Consider another feasible solution $f'_{ij}$. From
(\ref{positivef}) and (\ref{plot}) we have
		\begin{equation}
f'_{ij} (\alpha_i - \beta_j - c_{ij}) \ge 0 \qquad \forall i, j
		\end{equation}
and therefore, summing over $i$ and $j$ and using (\ref{rowsum})
and (\ref{columnsum}):
		\begin{equation}
\sum_i a_i \alpha_i - \sum_j b_j \beta_j 
- \sum_i \sum_j f'_{ij} c_{ij} \ge 0.
		\end{equation}
Comparing with (\ref{bidon1}), we have
		\begin{equation}
\sum_i \sum_j f'_{ij} c_{ij} \le \sum_i \sum_j F_{ij} c_{ij}
		\end{equation}
which shows that $F_{ij}$ is optimal. \endproof

Incidentally, (\ref{bidon1}) shows that the ``cost''
$\sum_i \sum_j F_{ij} c_{ij}$ of the optimal solution is equal
to the potential energy $U$ of the corresponding equilibrium.
\section{Numerical simulation}
\label{s:num-sim}

\subsection{Method}

The analog computer described in the previous Section could be
built in principle, but for large values of $m$ and $n$ this would
be impractical.  Instead, we will {\em simulate} on a digital
computer the behaviour of the analog computer, as it progressively
settles into an equilibrium.

It would be possible to simulate the evolution of the mechanical
system as described in the previous Section, i.e. to remove
suddenly the two stops $P$ and $Q$ and let the system evolve
freely until it has found an equilibrium. In that case, however,
all rods would interact more or less simultaneously, and the
simulation would be somewhat complex; essentially we would have to
solve an N-body problem. It turns out
to be simpler and also more efficient to guide the system through
a more controlled and orderly evolution. This is permitted
because, as shown by Theorem~\ref{th:equiv}, all we need in order
to solve the transportation problem is to find an equilibrium; how we
arrive at it is irrelevant.

Many algorithms can be imagined. Here we will only describe one of
the simplest methods, which was found to work well, although it is
not always the most efficient in terms of computing time (see
below Section~\ref{s:implementation}). We give here an informal
description of the algorithm, based on physical intuition; a
rigorous derivation will be presented in Section~\ref{s:formal}.

The stops $P$ and $Q$ are not removed. Instead, the stop $P$ is
held fixed during the whole process, and the stop $Q$ is slowly
lowered from its initial position.  The velocity of descent is
smaller than the minimal non-zero velocity of any aggregate,
$|M|_{\rm min}$, so that the evolution is fully controlled by the
motion of the stop $Q$. At any given time the system is in
quasi-equilibrium: if $Q$ stops then nothing moves anymore.

The evolution of the system will be studied in detail below.  It
ends when the whole system of rows and columns comes to rest. The
stop $Q$, continuing its descent, ceases then to be in contact
with any row and can be removed.  The force exerted by the stop
$P$ on any column still in contact with it is then zero, and that
stop can also be removed.  Thus an equibrilium has been reached,
 from which the optimal solution can be read.
\subsection{Graph representation}
\label{s:graph}

\begin{figure}[hbtp]
		$$\begin{picture}(200,51)


\put(5,10){\line(1,1){30}}
\put(35,10){\line(0,1){30}}
\put(65,10){\line(-1,1){30}}

\coinbr{14,21}{a_1}
\coinbl{37,21}{a_2}
\coinbl{56,21}{a_3}

\carre{0,0}{A_1}
\carre{30,0}{A_2}
\carre{60,0}{A_3}

\rond{35,46}{Q}


\put(105,11){\line(3,2){45}}
\put(135,11){\line(1,2){15}}
\put(165,11){\line(-1,2){15}}
\put(195,11){\line(-3,2){45}}

\coinbr{117,21}{b_1}
\coinbr{137,21}{b_2}
\coinbl{163,21}{b_3}
\coinbl{183,21}{b_4}

\carre{145,41}{P}

\rond{105,5}{B_1}
\rond{135,5}{B_2}
\rond{165,5}{B_3}
\rond{195,5}{B_4}

		\end{picture}$$
\caption{\label{f:initial}
Graph of the initial state.}
		\end{figure}
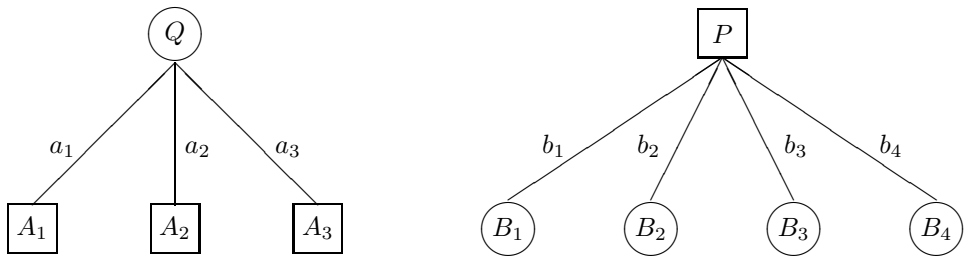

The state of the system at any given time can be conveniently
represented by a graph, as follows.  Each rod is a node of the
graph; we represent rows by squares and columns by circles.
An edge always joins a square to a circle: the graph is {\em
bipartite}. An edge is present between row $i$ and column $j$
when the row and the column are in contact through the stud
$c_{ij}$, i.e. when
		\begin{equation}
\alpha_i - \beta_j = c_{ij}.
					\label{contact}
\end{equation} 
The stops are also represented by nodes, and their contacts with
rods are similarly represented by edges. In order to preserve the
bipartite property, the stop $P$ must then be a square while the
stop $Q$ is a circle. As an illustration, Fig.~\ref{f:initial}
represents the initial state of the system, before any row has
been lowered. Note that the graph is a purely topological
representation: the position of a symbol in the graph has nothing
to do with the position of the corresponding rod in physical
space.

The force transmitted downwards from a row to a column, or from a
row to the $Q$ stop, or from the $P$ stop to a column, can be
written beside the corresponding edge. In the initial state, only
the last two kinds of forces are present; they have the values
indicated on Fig.~\ref{f:initial}. Note that this force is always
positive or zero.  At any given time in the procedure, each row is
in equilibrium; therefore the sum of the forces emanating from it
must equal the absolute value of its weight $a_i$. Similarly, each
column is in equilibrium, and the sum of the forces received by it
must equal its buoyancy $b_j$.

Can this graph have cycles ? A cycle will be an even sequence of
alternating rows and columns:
		\begin{equation}
i_1, j_1, i_2, j_2, \ldots, i_p, j_p.
		\end{equation}
 From (\ref{contact}) we have then
		\begin{equation}
c_{i_1j_1} - c_{i_2j_1} + c_{i_2j_2} - c_{i_3j_2} + \ldots 
+ c_{i_pj_p} - c_{i_1j_p} = 0.
					\label{cycle}
		\end{equation}
In order to simplify the exposition, we make the following
assumption (which will be removed in Section~\ref{s:formal}):
		\begin{assumption}
\label{as:c-incomm}
The $c_{ij}$ are such that there are no linear relations of the
form (\ref{cycle}) between them.
		\end{assumption}

Then the graph has no cycles and is a {\em forest}, or a
collection of trees. Each tree corresponds to an aggregate as
defined in Section~\ref{s:analog}.

It will be convenient to assume also, for the time being, that no
subset of rows and columns can be in equilibrium. This can be
expressed as:
			\begin{assumption}
\label{as:ab-incomm}
Let $I$ be a subset of $\{1, \dots, m\}$ and $J$ a
subset of $\{1, \dots, n\}$. The relation
		\begin{equation}
\sum_{i \in I} a_i = \sum_{j \in J} b_j
					\label{equal-subsets}
		\end{equation}
is true only in two cases: (i) $I = J = \emptyset$; (ii) $I = \{1,
\dots, m\}$ and $J = \{1, \dots, n\}$.
		\end{assumption}

Note that this excludes in particular $a_i = 0$ or $b_j = 0$:
rows and columns must have positive weights and buoyancies.

We remark that these two assumptions are satisfied in principle in
the generic case, when the $a_i$, $b_j$, $c_{ij}$ are real numbers
with arbitrary values. In practice, however, these numbers are
often integers with a restricted range and the assumptions are
frequently violated.

\begin{figure}[hbtp]
		$$\begin{picture}(170,144)(10,10)


\put(15,112){\line(3,4){15}}
\put(45,112){\line(-3,4){15}}
\put(30,82){\line(3,4){15}}
\put(60,82){\line(-3,4){15}}
\put(15,50){\line(3,4){15}}
\put(45,50){\line(-3,4){15}}

\carre{10,102}{A_6}
\carre{40,102}{A_5}
\carre{10,40}{A_2}
\carre{40,40}{A_3}

\rond{30,138}{Q}
\rond{30,76}{B_4}
\rond{60,76}{B_6}


\put(105,113){\line(3,2){30}}
\put(135,113){\line(0,1){20}}
\put(165,113){\line(-3,2){30}}
\put(135,81){\line(0,1){20}}
\put(120,51){\line(3,4){15}}
\put(150,51){\line(-3,4){15}}
\put(120,19){\line(0,1){20}}

\carre{130,133}{P}
\carre{130,71}{A_1}
\carre{115,9}{A_4}

\rond{105,107}{B_2}
\rond{135,107}{B_3}
\rond{165,107}{B_7}
\rond{120,45}{B_1}
\rond{150,45}{B_5}

		\end{picture}$$
\caption{\label{f:twotrees}
Example of a graph. Left: moving tree. Right: fixed tree.}
		\end{figure}
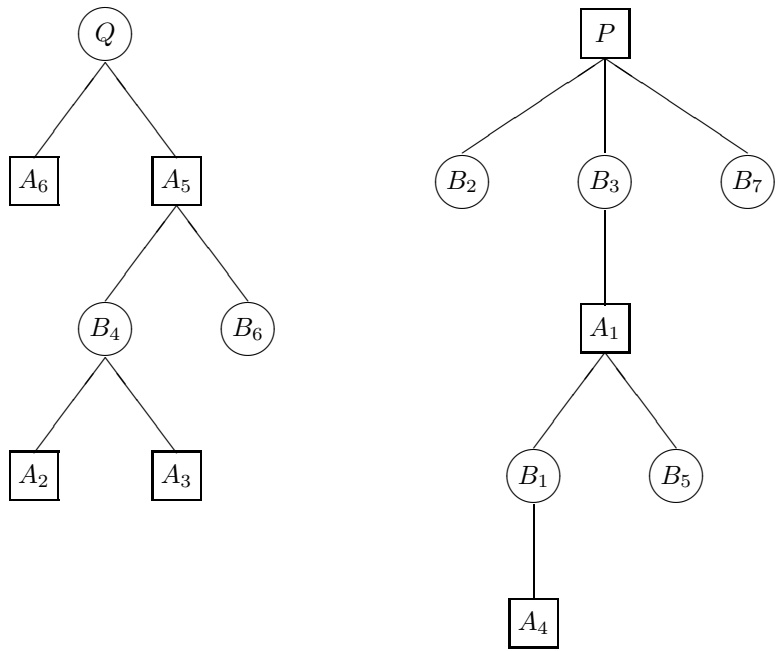

It will be proved in Section~\ref{s:formal} that the graph always
consists of two trees, except at particular instants of time where
they fuse into a single tree (see below
Section~\ref{s:cContact-and-rearrangement}). One of them contains
$Q$ and will be called {\em moving tree}. The other contains $P$
and will be called {\em fixed tree}.  We will represent the two
trees with the usual hierarchical representation of trees
(\cite{Knu73a}, Section~2.3), taking the stop as root for each
tree.  Fig.~\ref{f:twotrees} shows an example.


It will be convenient to extend the usual terminology
of parents and children by specifying that {\em rows are male}
and {\em columns are female}. The stop $Q$ is female and the stop
$P$ is male. To recapitulate:
		\begin{eqnarray}
\left. \begin{array}{c} \mbox{row} \\ \mbox{stop $P$}
\end{array} \right\} 
& = & \begin{picture}(15,15) \carre{0,-4}{} \end{picture}
= \mbox{male}, \nonumber \\
\left. \begin{array}{c} \mbox{column} \\ \mbox{stop $Q$}
\end{array} \right\} 
& = & \begin{picture}(15,15) \rond{6,2}{} \end{picture}
= \mbox{female}. 
		\end{eqnarray}

In each tree, sex is thus alternating from one generation
to the next.  A row has one mother and any number of daughters,
while a column has one father and any number of sons.  The stops
themselves have no parents.
\subsection{Contact and rearrangement}
\label{s:cContact-and-rearrangement}

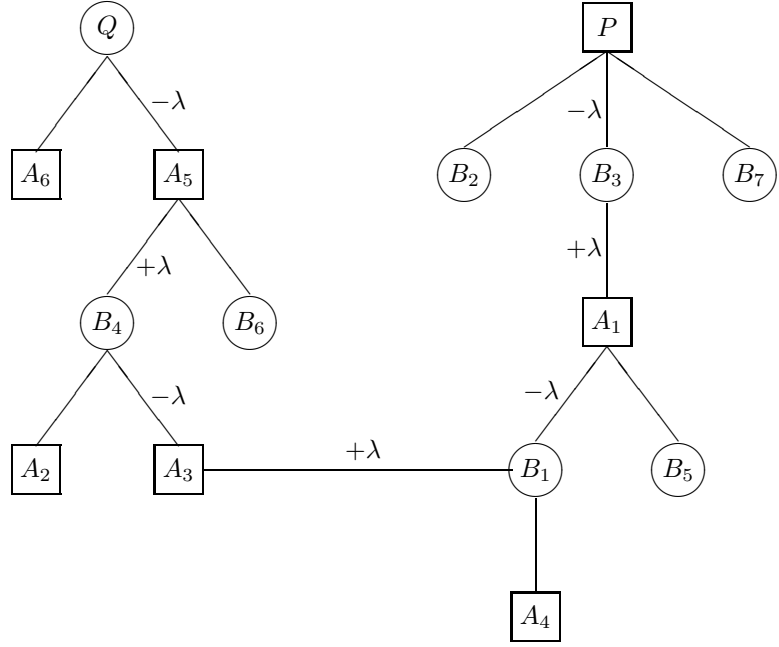
\begin{figure}[hbtp]
		$$\begin{picture}(170,140)(10,10)


\put(15,112){\line(3,4){15}}
\put(45,112){\line(-3,4){15}}
\put(30,82){\line(3,4){15}}
\put(60,82){\line(-3,4){15}}
\put(15,50){\line(3,4){15}}
\put(45,50){\line(-3,4){15}}

\carre{10,102}{A_6}
\carre{40,102}{A_5}
\carre{10,40}{A_2}
\carre{40,40}{A_3}

\rond{30,138}{Q}
\rond{30,76}{B_4}
\rond{60,76}{B_6}


\put(105,113){\line(3,2){30}}
\put(135,113){\line(0,1){20}}
\put(165,113){\line(-3,2){30}}
\put(135,81){\line(0,1){20}}
\put(120,51){\line(3,4){15}}
\put(150,51){\line(-3,4){15}}
\put(120,19){\line(0,1){20}}

\carre{130,133}{P}
\carre{130,71}{A_1}
\carre{115,9}{A_4}

\rond{105,107}{B_2}
\rond{135,107}{B_3}
\rond{165,107}{B_7}
\rond{120,45}{B_1}
\rond{150,45}{B_5}


\put(50,45){\line(1,0){65}}

\coinbl{39,121}{-\lambda}
\cointl{36,90}{+\lambda}
\coinbl{39,59}{-\lambda}
\coinbl{80,47.5}{+\lambda}
\coinbr{125,60}{-\lambda}
\coinbr{134,90}{+\lambda}
\coinbr{134,119}{-\lambda}

		\end{picture}$$
\caption{\label{f:contact}
Contact between the moving and the fixed tree, and readjustment of
the forces. Here the numbers $+\lambda$ and $-\lambda$ represent
the changes in the forces, rather than the forces themselves.}
		\end{figure}

As the stop $Q$ goes down, the rows and columns which belong
to the moving tree move with it. This continues until a contact
is made between the two trees.  Since rows always remain above
columns in physical space, this contact happens necessarily
between a moving row $A_{i_c}$ and a fixed column $B_{j_c}$. A
new edge is created and the two trees are temporarily fused into a
single tree. A single contact is made, because two simultaneous
contacts would imply a cycle, in contradiction to Assumption
\ref{as:c-incomm}.  Fig.~\ref{f:contact} shows an example,
posterior in time to Fig.~\ref{f:twotrees},  where contact is made
between the row $A_3$ and the column $B_1$. (Note: in the
example of Fig.~\ref{f:contact}, the row and the column which
make contact happen to be at the same level in their respective
trees; this need not be so in general.)


As a result of the contact, the forces change.  To understand what
happens, it is convenient to imagine that the studs are slightly
elastic, so that the change does not happen all at once, but
progressively over a small interval of time.  The force $\lambda$
along the newly created edge, which is initially zero, increases
as the moving row moves down. This induces other changes in
neighbouring edges.  We consider the path from $Q$ to $P$ made by
(i) the path from $Q$ to $A_{i_c}$ in the moving tree; (ii) the
newly created edge from $A_{i_c}$ to $B_{j_c}$; (iii) the path
 from $B_{j_c}$ to $P$ in the fixed tree. We call this the {\em
main path}. In Fig.~\ref{f:contact}, for instance, the main path
is $Q A_5 B_4 A_3 B_1 A_1 B_3 P$. The graph can then be viewed as
made up of the main path, plus a number of lateral branches. The
lateral branches are not involved in the readjustment of the
forces; each of them is attached to the main path by a single edge
and the force along that edge equals the weight or the buoyancy of
the branch, which does not change.  Therefore, only forces along
the main path can change. Each node must remain in equilibrium;
therefore all changes have the same modulus $\lambda$ and
alternate in sign along the main path, as shown by
Fig.~\ref{f:contact}. (Note in particular that the forces of
contact with the two stops decrease).

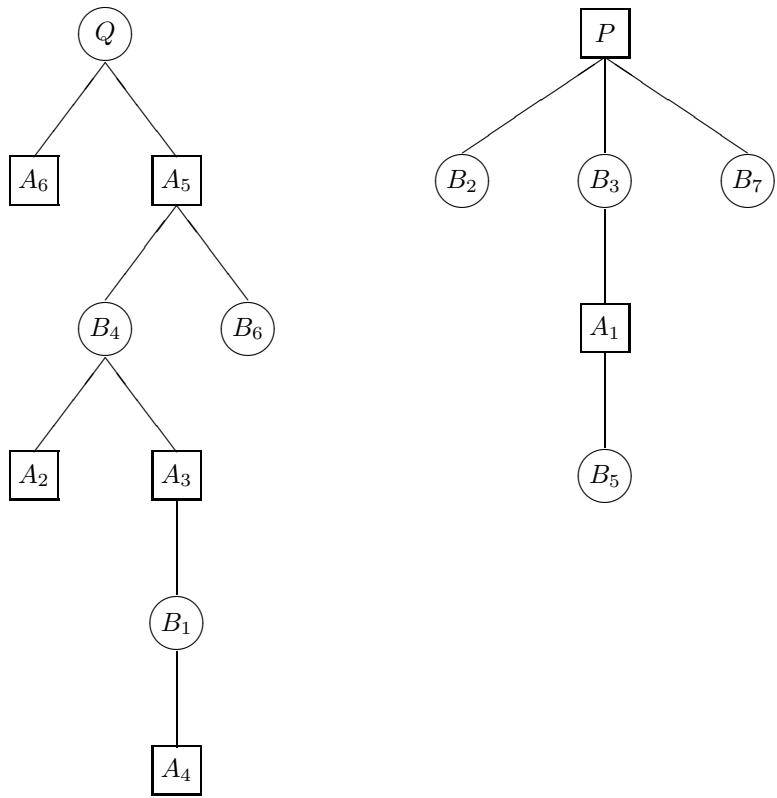
\begin{figure}[hbtp]
		$$\begin{picture}(170,174)(10,10)


\put(15,142){\line(3,4){15}}
\put(45,142){\line(-3,4){15}}
\put(30,112){\line(3,4){15}}
\put(60,112){\line(-3,4){15}}
\put(15,80){\line(3,4){15}}
\put(45,80){\line(-3,4){15}}
\put(45,50){\line(0,1){20}}
\put(45,18){\line(0,1){20}}

\carre{10,132}{A_6}
\carre{40,132}{A_5}
\carre{10,70}{A_2}
\carre{40,70}{A_3}
\carre{40,8}{A_4}

\rond{30,168}{Q}
\rond{30,106}{B_4}
\rond{60,106}{B_6}
\rond{45,44}{B_1}


\put(105,143){\line(3,2){30}}
\put(135,143){\line(0,1){20}}
\put(165,143){\line(-3,2){30}}
\put(135,111){\line(0,1){20}}
\put(135,81){\line(0,1){20}}

\carre{130,163}{P}
\carre{130,101}{A_1}

\rond{105,137}{B_2}
\rond{135,137}{B_3}
\rond{165,137}{B_7}
\rond{135,75}{B_5}

		\end{picture}$$
\caption{\label{f:capture}
Example of a capture of a part of the fixed tree by the moving
tree.}
		\end{figure}

This continues until one of the decreasing forces becomes zero.
(Two forces cannot vanish simultaneoulsy, because the intermediate
tree would have zero weight, in contradiction to
Assumption~\ref{as:ab-incomm}).  The chain breaks then at the
corresponding edge, which disappears, and we have again two
separate trees.  Thus, the whole episode ends in a {\em capture}
of a part of one tree by the other. The capture can occur in
either direction, depending on where is the weakest link of the
main path. For instance if the weakest link in
Fig.~\ref{f:contact} is between $B_1$ and $A_1$, the branch of the
fixed tree with head $B_1$ is captured by the moving tree and we
obtain Fig.~\ref{f:capture}. $B_1$ ceases to be in contact with
$A_1$ as the moving tree continues its descent.


The new moving tree continues to go down with the stop $Q$.
Eventually a new contact is made, and one of the trees captures a
part of the other.  This goes on until the moving tree is reduced
to the stop $Q$ alone. It can be shown that this always happens
after a finite number of captures (see
Appendix~\ref{s:number-of-cycles}). Only the fixed tree remains,
now containing all rows and all columns, and we have the sought
equilibrium.
\clearpage
\section{Example}
\label{s:example}


\setlength{\unitlength}{1.2pt}
\renewcommand{\coinbl}[2]{\put(#1){\makebox(0,0)[bl]%
{\small $#2$}}}
\renewcommand{\coinbr}[2]{\put(#1){\makebox(0,0)[br]%
{\small $#2$}}}
\renewcommand{\cointl}[2]{\put(#1){\makebox(0,0)[tl]%
{\small $#2$}}}
\renewcommand{\cointr}[2]{\put(#1){\makebox(0,0)[tr]%
{\small $#2$}}}

We exhibit here the step-by-step progress of the algorithm on a
simple example with $m = 3$, $n = 4$. Table~\ref{t:ex-param} shows
the values of the given coefficients $a_i$, $b_j$, and $c_{ij}$.
Note that the condition (\ref{zerosum}) is verified.

		\begin{table}[h]
\caption{Values of the parameters $a_i$ (left column), $b_j$
(top row), and $c_{ij}$ for the example problem.}
\label{t:ex-param}
		$$\begin{tabular}{|c|cccc|}
\hline
   & 44 & 52 & 13 & 37 \\
\hline
86 & 26 & 64 & 33 & 62 \\
 4 & 63 & 27 & 13 & 14 \\
56 & 94 &  4 &  4 & 52 \\
\hline
		\end{tabular}$$
		\end{table}

\begin{figure}[p]
a\hspace{10mm}%
		\begin{tabular}{|cccc|}
\hline
74 & 36 & 67 & 38 \\
37 & 73 & 87 & 86 \\
 6 & 96 & 96 & 48 \\
\hline
		\end{tabular}

\vspace{5mm}
\raisebox{9mm}{b \hspace{17mm} $d = 6$}
\hfill
		\begin{picture}(206,44.5)(-60,196.5)
\rond{-25,235}{Q}
\put(-55,206.5){\line(4,3){30}}
\put(-25,206.5){\line(0,1){22.5}}
\put(5,206.5){\line(-4,3){30}}
\coinbr{-50,212.5}{86}
\coinbl{-23.5,211.5}{4}
\coinbl{0,212.5}{56}
\carre{-60,196.5}{1}
\carre{-30,196.5}{2}
\carre{0,196.5}{3}

\carre{90,230}{P}
\put(50,207.5){\line(2,1){45}}
\put(80,207.5){\line(2,3){15}}
\put(110,207.5){\line(-2,3){15}}
\put(140,207.5){\line(-2,1){45}}
\coinbr{58,214}{44}
\coinbr{82,212}{52}
\coinbl{108,212}{13}
\coinbl{132,214}{37}
\rond{50,201.5}{1}
\rond{80,201.5}{2}
\rond{110,201.5}{3}
\rond{140,201.5}{4}
		\end{picture}
\vspace{5mm}

c\hspace{10mm}%
		\begin{tabular}{|cccc|}
\hline
68 & 30 & 61 & 32 \\
31 & 67 & 81 & 80 \\
 0 & 90 & 90 & 42 \\
\hline
		\end{tabular}

\raisebox{15mm}{d \hspace{17mm} $d = 30$}
\hfill
		\begin{picture}(206,74)(-60,167)
\rond{-25,235}{Q}
\put(-55,209){\line(3,2){30}}
\put(-25,209){\line(0,1){20}}
\put(5,209){\line(-3,2){30}}
\coinbr{-50,215}{86}
\coinbl{-23.5,214}{4}
\coinbl{0,215}{12}
\carre{-60,199}{1}
\carre{-30,199}{2}
\carre{0,199}{3}
\put(5,179){\line(0,1){20}}
\coinbl{7,187}{44}
\rond{5,173}{1}

\carre{90,230}{P}
\put(65,210){\line(3,2){30}}
\put(95,210){\line(0,1){20}}
\put(125,210){\line(-3,2){30}}
\coinbr{70,216}{52}
\coinbl{96.5,215}{13}
\coinbl{120,216}{37}
\rond{65,204}{2}
\rond{95,204}{3}
\rond{125,204}{4}
		\end{picture}

e\hspace{10mm}%
		\begin{tabular}{|cccc|}
\hline
68 &  0 & 31 &  2 \\
31 & 37 & 51 & 50 \\
 0 & 60 & 60 & 12 \\
\hline
		\end{tabular}

\raisebox{15mm}{f \hspace{17mm} $d = 2$}
\hfill
		\begin{picture}(206,74)(-60,167)
\rond{-25,235}{Q}
\put(-55,209){\line(3,2){30}}
\put(-25,209){\line(0,1){20}}
\put(5,209){\line(-3,2){30}}
\coinbr{-50,215}{34}
\coinbl{-23.5,214}{4}
\coinbl{0,215}{12}
\carre{-60,199}{1}
\carre{-30,199}{2}
\carre{0,199}{3}
\put(-55,179){\line(0,1){20}}
\put(5,179){\line(0,1){20}}
\coinbl{-53,187}{52}
\coinbl{7,187}{44}
\rond{-55,173}{2}
\rond{5,173}{1}

\carre{90,230}{P}
\put(80,210){\line(3,4){15}}
\put(110,210){\line(-3,4){15}}
\coinbr{83,216}{13}
\coinbl{107,216}{37}
\rond{80,204}{3}
\rond{110,204}{4}
		\end{picture}

\caption{\label{f:evol} Exemple. Left: the distances 
$\gamma_{ij} = \alpha_i - \beta_j - c_{ij}$ between the lines and
the studs. Right: moving tree and fixed tree.}
		\end{figure}
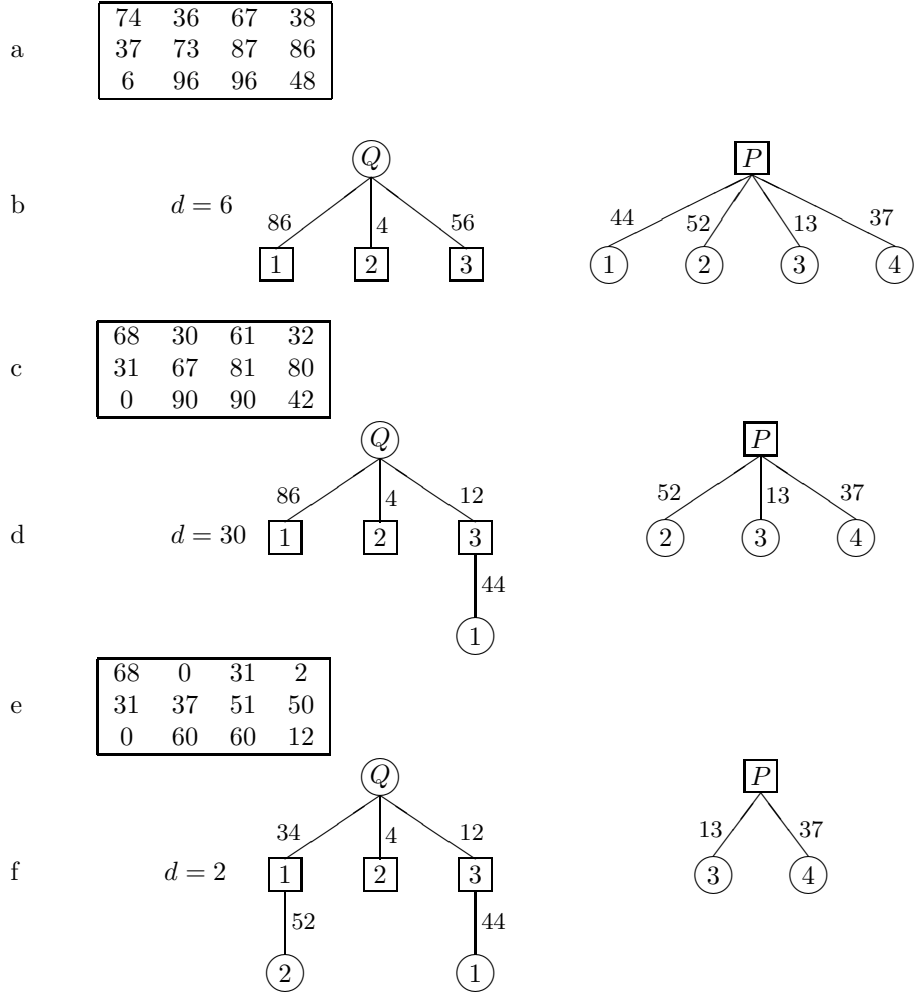

When executing the algorithm by hand, it is convenient to keep
track of the distances between the rows and the studs, i.e. the
quantities 
\begin{equation}
\gamma_{ij} = \alpha_i - \beta_j - c_{ij}. 
				\label{defgamma}
\end{equation}
One can then easily determine where the next contact will take
place.  The distances $\gamma_{ij}$ are shown on the left in
Fig.~\ref{f:evol}, while the graph (moving tree and fixed tree) is
shown on the right. Only the indices $i$ or $j$ of the rows and
columns are indicated; the type is indicated by the symbol (square
for a row, circle for a column). Evolution proceeds from top to
bottom; successive steps are represented in lines labelled a, b,
c, \dots. Lines~b, d, f, \dots, correspond to a descent of the
moving tree; the distances change, while the forces and the trees
remain fixed. Conversely, lines c, e, g, \dots, correspond to a
readjustment of the forces and of the trees, during which the
distances do not change.

Initially we set the height of the rows and columns at $\alpha_i =
100$ and $\beta_j = 0$. The corresponding distances $\gamma_{ij}$
are then obtained by complementing to 100 the values of
table~\ref{t:ex-param} and are shown in Fig.~\ref{f:evol}, line~a.
The initial moving and fixed trees are set up as indicated in
Section~\ref{s:graph}, Fig.~\ref{f:initial}, and are shown in
line~b.

All rows are moving and all columns are fixed, therefore all
distances $\gamma_{ij}$ decrease. From line~a, we immediately find
that the moving tree can descend a distance $d = 6$; a contact is
then made between row 3 and column 1. The new distances are shown
in line~c.

We now readjust the forces and the trees. The main path is: stop
$Q$ $-$ row~3 $-$ column~1 $-$ stop $P$. The weakest link is
between the column 1 and the stop $P$, with a force 44. Therefore
the column 1 is captured by the moving tree. The forces along the
main path change by $\lambda = \pm 44$.  The new trees and the new
forces are represented in line~d.

All rows are still moving; in addition, column 1 is also moving.
Therefore the distances $\gamma_{ij}$ remain fixed for $j = 1$
(first column of matrix) and decrease for $j \in \{2, 3, 4\}$.
 From line~c we find then that the distance of descent is $d = 30$.
Contact is made between row 1 and column 2. The new distances are
shown in line~e.

Column 2 is captured by the moving tree and we obtain line~f. Now
the distances decrease for $j \in \{3, 4\}$. After a descent of $d
= 2$, contact is made between row 1 and column 4. The new
distances are shown in line~g.

This time the weakest link is between $Q$ and row 1.  Therefore
row 1, and its daughter the column 2, are captured by the fixed
tree. The new trees are represented in line~h.

Now only the distances $\gamma_{ij}$ with $i \in \{2,3\}$ and $j
\in \{2,3,4\}$ are decreasing. Therefore we have a descent of $d =
10$. Note that $\gamma_{11}$ is {\em increasing}, since row 1 is
fixed and column 1 is moving. The other distances remain fixed.
The new distances are shown in line~i. Contact is made between
row~3 and column~4.

The evolution continues. In line~i, the column 4 and its two
descendants are captured by the moving tree. In line~k, we observe
a more complex event, involving the capture of a large piece of
the main path and a drastic reorganization of the trees. Finally,
in line~m the last remnant of the moving tree is captured. 

We reach line~n, where only the fixed tree remains. The forces
$f_{ij}$ along the edges of that tree give the solution of the
transportation problem. They can be rewritten in matrix form
(Table~\ref{t:ex-sol}).
\clearpage
\clearpage

g\hspace{10mm}%
		\begin{tabular}{|cccc|}
\hline
68 &  0 & 29 &  0 \\
31 & 37 & 49 & 48 \\
 0 & 60 & 58 & 10 \\
\hline
		\end{tabular}

\raisebox{21mm}{h \hspace{17mm} $d = 10$}
\hfill
		\begin{picture}(161,105)(-30,136)
\rond{-10,235}{Q}
\put(-25,209){\line(3,4){15}}
\put(5,209){\line(-3,4){15}}
\coinbr{-22,215}{4}
\coinbl{2,215}{12}
\carre{-30,199}{2}
\carre{0,199}{3}
\put(5,179){\line(0,1){20}}
\coinbl{7,186}{44}
\rond{5,173}{1}

\carre{90,230}{P}
\put(80,210){\line(3,4){15}}
\put(110,210){\line(-3,4){15}}
\coinbr{83,216}{13}
\coinbl{107,216}{3}
\rond{80,204}{3}
\rond{110,204}{4}
\put(110,178){\line(0,1){20}}
\coinbl{112,186}{34}
\carre{105,168}{1}
\put(110,148){\line(0,1){20}}
\coinbl{112,156}{52}
\rond{110,142}{2}
		\end{picture}

i\hspace{10mm}%
		\begin{tabular}{|cccc|}
\hline
78 &  0 & 29 &  0 \\
31 & 27 & 39 & 38 \\
 0 & 50 & 48 &  0 \\
\hline
		\end{tabular}

\raisebox{27mm}{j \hspace{17mm} $d = 29$}
\hfill
		\begin{picture}(176,136)(-45,105)
\rond{-10,235}{Q}
\put(-25,209){\line(3,4){15}}
\put(5,209){\line(-3,4){15}}
\coinbr{-22,215}{4}
\coinbl{2,215}{9}
\carre{-30,199}{2}
\carre{0,199}{3}
\put(-10,179){\line(3,4){15}}
\put(20,179){\line(-3,4){15}}
\coinbr{-7,185}{44}
\coinbl{17,185}{3}
\rond{-10,173}{1}
\rond{20,173}{4}
\put(20,147){\line(0,1){20}}
\coinbl{22,155}{34}
\carre{15,137}{1}
\put(20,117){\line(0,1){20}}
\coinbl{22,125}{52}
\rond{20,111}{2}

\carre{90,230}{P}
\put(95,210){\line(0,1){20}}
\coinbl{97,218}{13}
\rond{95,204}{3}
		\end{picture}

k\hspace{10mm}%
		\begin{tabular}{|cccc|}
\hline
78 &  0 &  0 &  0 \\
31 & 27 & 10 & 38 \\
 0 & 50 & 19 &  0 \\
\hline
		\end{tabular}

\vspace{10mm}
\centerline{Figure~\thefigure\ (continued).}
\newpage

\raisebox{34mm}{l \hspace{17mm} $d = 10$}
\hfill
		\begin{picture}(176,167)(-45,74)
\rond{-10,235}{Q}
\put(-10,209){\line(0,1){20}}
\coinbl{-8,217}{4}
\carre{-15,199}{2}

\carre{90,230}{P}
\put(95,210){\line(0,1){20}}
\coinbl{97,218}{4}
\rond{95,204}{3}
\put(95,178){\line(0,1){20}}
\coinbl{97,186}{9}
\carre{90,168}{1}
\put(80,148){\line(3,4){15}}

\put(110,148){\line(-3,4){15}}
\coinbr{83,154}{25}
\coinbl{107,154}{52}
\rond{80,142}{4}
\rond{110,142}{2}
\put(80,116){\line(0,1){20}}
\coinbl{82,124}{12}
\carre{75,106}{3}
\put(80,86){\line(0,1){20}}
\coinbl{82,94}{44}
\rond{80,80}{1}
		\end{picture}

m\hspace{10mm}%
		\begin{tabular}{|cccc|}
\hline
78 &  0 &  0 &  0 \\
21 & 17 &  0 & 28 \\
 0 & 50 & 19 &  0 \\
\hline
		\end{tabular}

\raisebox{34mm}{n \hspace{15mm} equilibrium}
\hfill
		\begin{picture}(176,167)(-45,74)
\carre{90,230}{P}
\put(95,210){\line(0,1){20}}
\coinbl{97,218}{0}
\rond{95,204}{3}
\put(80,178){\line(3,4){15}}
\put(110,178){\line(-3,4){15}}
\coinbr{83,184}{4}
\coinbl{107,184}{9}
\carre{75,168}{2}
\carre{105,168}{1}
\put(95,148){\line(3,4){15}}
\put(125,148){\line(-3,4){15}}
\coinbr{98,154}{25}
\coinbl{122,154}{52}
\rond{95,142}{4}
\rond{125,142}{2}
\put(95,116){\line(0,1){20}}
\coinbl{97,124}{12}
\carre{90,106}{3}
\put(95,86){\line(0,1){20}}
\coinbl{97,94}{44}
\rond{95,80}{1}
		\end{picture}

\vspace{10mm}
\centerline{Figure~\thefigure\ (continued).}

\setlength{\unitlength}{1.8pt}
\renewcommand{\coinbl}[2]{\put(#1){\makebox(0,0)[bl]{$#2$}}}
\renewcommand{\coinbr}[2]{\put(#1){\makebox(0,0)[br]{$#2$}}}
\renewcommand{\cointl}[2]{\put(#1){\makebox(0,0)[tl]{$#2$}}}
\renewcommand{\cointr}[2]{\put(#1){\makebox(0,0)[tr]{$#2$}}}
\newpage
		\begin{table}[htb]
\caption{Solution of the transportation problem defined by
Table~\protect\ref{t:ex-param}.}
\label{t:ex-sol}
		$$\begin{tabular}{|cccc|}
\hline
 0 & 52 &  9 & 25 \\
 0 &  0 &  4 &  0 \\
44 &  0 &  0 & 12 \\
\hline
		\end{tabular}$$
		\end{table}
\section{Formal definition and justification of the algorithm}
\label{s:formal}

\subsection{Definitions}

In the present Section, we give a more rigorous definition of the
numerical algorithm, and we prove that it solves the transportation
problem. We also drop the restrictions introduced in
Section~\ref{s:num-sim}: assumptions \ref{as:c-incomm} and
\ref{as:ab-incomm} do not have to be satisfied any more, i.e.
relations of the form (\ref{cycle}) and (\ref{equal-subsets}) are
allowed.

This Section is independent of the description of the mechanical
model in Section~\ref{s:analog}, and also of the informal
description of the numerical algorithm in Section~\ref{s:num-sim}.
It will be sometimes convenient to use names which are reminiscent
of the origin of the algorithm, such as ``descent'' or ``force'';
the reasoning, however, will be purely mathematical.

In what follows, unless otherwise specified, $i$ is always
understood to take all values from 1 to $m$ and $j$ to take all
values from 1 to $n$.

We want to solve a transportation problem defined by given $m$, $n$,
$a_i$, $b_j$, $c_{ij}$ (Section~\ref{s:problem}). The algorithm
operates on the following collection of objects:

		\begin{itemize}
\item A graph with $m + n + 2$
nodes labelled $A_1, \dots, A_m, B_1, \dots, B_n, P, Q$.  

This set of nodes remains invariant during the course of the
computation. On the other hand, the set of edges varies.

\item To each $A_i$ node is associated a variable
number $\alpha_i$. To each $B_j$ node is similarly associated a
variable number $\beta_j$.

\item To each edge is also associated a variable number, which
will be called a {\em force}.

		\end{itemize}

\subsection{Properties}

In the following Section, we will define the operation of the
algorithm; simultaneously, we will prove that the following
properties hold throughout the computation.

		\begin{description}

\item [P1] Only the following kinds of edges are allowed: between
an $A_i$ and a $B_j$, between an $A_i$ and $Q$, and between $P$
and a $B_j$.

It follows that the graph is {\em bipartite}, the two subsets of
nodes being $\{A_1, \dots, A_m, P\}$ and $\{B_1, \dots, B_n,
Q\}$.

\item [P2] The graph is a forest, consisting of one or two trees.

\item [P3] When there are two trees, one of them
contains $P$ and at least one other node. The other tree contains
$Q$ and at least one other node. They will be called
respectively $T_f$ or {\em fixed tree} and $T_m$ or {\em moving
tree}.

\item [P4] $\gamma_{ij} \ge 0$. ($\gamma_{ij}$ is the {\sl
distance} defined by (\ref{defgamma})).

\item [P5] If there is an edge
between $A_i$ and $B_j$, then $\gamma_{ij} = 0$.

\item [P6] Forces are positive or zero.

\item [P7] The sum of the forces on the edges
adjacent to node $A_i$ equals $a_i$. The sum of the forces on the
edges adjacent to node $B_j$ equals $b_j$.

		\end{description}

\subsection{Algorithm}
\label{s:algorithm}

The algorithm consists in a succession of steps, described in the
following Sections. Step 1 is executed only once. Then a main
loop, made of steps 2 to 5, is executed a number of times; each
execution will be called a {\em cycle}.
\subsubsection*{Step 1: Initialize}

We set up the initial state of the graph and of the associated
values as follows. An edge is established between $Q$ and each of
the $A_i$, with associated force $a_i$, and between $P$ and each
of the $B_j$, with associated force $b_j$ (see
Fig.~\ref{f:initial}).  We set $\alpha_i = c_{\rm sup}$, with
$c_{\rm sup}$ satisfying (\ref{defcsup}), and $\beta_j = 0$. It is
easily verified that properties {\bf P1} to {\bf P7} hold. There
are two trees.
\subsubsection*{Step 2: Descent}

Since there are two trees, property {\bf P3} ensures that there is
at least one node in the moving tree other than $Q$.  Since $Q$
can be connected only to $A_i$ nodes, the moving tree includes at
least one $A_i$ node. Similarly, the fixed tree includes at least
one $B_j$ node. Therefore the following minimum exists and can be
computed: 		\begin{equation} d = \min_{A_i \in T_m,
B_j \in T_f} \gamma_{ij}.
\label{d-def} 		\end{equation} From property {\bf P4} we
have: $d \ge 0$.

Next we effect the ``descent of the moving tree'':
		\begin{eqnarray}
\alpha_i & := & \alpha_i - d \qquad {\rm for\ all\ } A_i \in T_m,
  \nonumber \\
\beta_j & := & \beta_j - d \qquad {\rm for\ all\ } B_j \in T_m.
					\label{descent} 
		\end{eqnarray}

Note that $d$ may be zero, in which case nothing changes.

We verify now that the properties still hold. Only the $\alpha_i$
and $\beta_j$ have changed, therefore we have only to examine
properties {\bf P4} and {\bf P5}. We consider first {\bf P4}. If
$\alpha_i$ and $\beta_j$ belong to the same tree, $\gamma_{ij}$
does not change. If $A_i \in T_f$ and $B_j \in T_m$, $\gamma_{ij}$
increases. Finally, if $A_i \in T_m$ and $B_j \in T_d$,
$\gamma_{ij}$ decreases by $d$, but remains positive or zero as a
consequence of (\ref{d-def}).

We verify also {\bf P5}: if there is an edge between $A_i$ and
$B_j$, these nodes belong to the same tree, and therefore
$\gamma_{ij}$ does not change.
\subsubsection*{Step 3: Contact}

We consider the pair of values $i = i_c$, $j - j_c$ which realized
the minimum $d$ in step 2, i.e. which were such that $A_{i_c} \in
T_m$, $B_{j_c} \in T_f$, and $\gamma_{i_c j_c} = d$ before the
descent. (If more than one pair (i, j) realized the minimum, we
select one of them arbitrarily).  From (\ref{descent}) we find
that there is now, after the descent: $\gamma_{i_c j_c} = 0$.

We add one edge between nodes $A_{i_c}$ and $B_{j_c}$, and we set
the associated force equal to zero.

We consider the properties. {\bf P1} is still satisfied since the
new edge is between an $A_i$ and a $B_j$ node.  Concerning {\bf
P2}, since we have linked one nodes of $T_m$ with one node of
$T_f$, the graph now consists of a single tree. {\bf P3} does not
apply any more. {\bf P4} is not affected by a change in the graph.
{\bf P5} is satisfied for the new edge since $\gamma_{i_c j_c} =
0$. {\bf P6} is satisfied since the new force is zero. Finally,
{\bf P7} still holds for nodes $A_{i_c}$ and$B_{j_c}$, again
because the new force is zero.
\subsubsection*{Step 4: Readjustment}

We define the {\em main path} as the oriented path from $Q$ to
$P$. This path is unique since the graph consists of a single
tree. The main path is made of three parts: (i) the path from $Q$
to $A_{i_c}$ in the previous moving tree; (ii) the newly created
edge from $A_{i_c}$ to $B_{j_c}$; (iii) the path from $B_{j_c}$ to
$P$ in the previous fixed tree. From property {\bf P1}, we deduce
that the first part has an odd number of edges (the graph is
bipartite, and $Q$ and $A_i$ belong to different subsets).
Similarly, the last part has an odd number of edges.  Thus, the
main path as a whole has an odd number of edges, which is at least
equal to 3. We number the edges along the main path, from $Q$ to
$P$, starting from 1. We note that the two end edges, adjacent to
$Q$ and $P$, are odd-numbered, and that the newly created edge is
even-numbered.

We compute the minimum $\lambda$ of the forces associated with the
odd-numbered edges. There is $\lambda \ge 0$ by virtue of
property {\bf P6}.

We readjust the forces along the main path, by adding $\lambda$ to
the forces associated with even-numbered edges and subtracting
$\lambda$ from the forces associated with odd-numbered edges.

Only the forces have changed, therefore we have only to examine
properties {\bf P6} and {\bf P7}. Property {\bf P6} is obviously
still true.  For every node $A_i$ or $B_j$ along the main path,
the changes of the forces associated with the two adjacent edges
on the main path cancel each other, so that {\bf P7} remains
true.
\subsubsection*{Step 5: Breaking}

We consider the odd-numbered edge of the main path which realized
the minimum in the previous step. (If more than one edge realized
the minimum, we select one of them arbitrarily). We will call it
the {\em breaking edge}. The force associated with that edge is
now zero.

We delete the breaking edge. This completes one cycle of the
algorithm.

Properties {\bf P1}, {\bf P5}, {\bf P6} are not affected since we
have simply removed an edge. The graph consists now again of two
trees: {\bf P2} is satisfied. Property {\bf P4} is not affected by
a change in the graph. Property {\bf P7} still holds for the two
nodes adjacent to the deleted edge since the force associated
with that edge was zero.

There remains to consider property {\bf P3}. Since the deleted
edge was on a path from $Q$ to $P$, it is still true that one of
the new trees contains $Q$ and the other contains $P$.  However it
can happen that $Q$ or $P$ is now an isolated node. We distinguish
two cases.

		\begin{enumerate}
\item The moving tree contains nodes other than $Q$, and the fixed
tree contains nodes other than $P$. Property {\bf P3} is
satisfied. We go back to step 2 for a new cycle.

\item \label{case:end} The moving tree contains $Q$ alone, or the
fixed tree contains $P$ alone.  This signals the end of the
computation.

		\end{enumerate}

It can be shown that case \ref{case:end} necessarily happens after
a finite number of steps: a definite upper bound on the number of
cycles can be computed (see Appendix~\ref{s:number-of-cycles}).
Thus the algorithm always terminates.
\subsection{Solution}

We show now that a solution of the transportation problem has been
obtained.  This is similar to the proof given at the end of
Section~\ref{s:analog}; the present derivation, however, makes no
reference to the mechanical model.

We will describe only the case where the moving tree contains $Q$
alone; the other case, where the fixed tree contains $P$ alone, is
treated in the same way, exchanging rows and columns.  We
consider the fixed tree. It contains all nodes $A_i$ and $B_j$, in
addition to $P$. We define $f_{ij}$ as follows: if there is an
edge between $A_i$ and $B_j$, $f_{ij}$ is equal to the associated
force; otherwise $f_{ij} = 0$. We also define $f_{*j}$ as follows:
if there is an edge between $P$ and $B_j$, $f_{*j}$ is equal to
the associated force; otherwise $f_{*j} = 0$.

By virtue of properties {\bf P1} and {\bf P7}, we have
		\begin{eqnarray}
\sum_j f_{ij} & = & a_i, \qquad (i = 1, \dots, m),
					\label{rowsum2} \\
\sum_i f_{ij} + f_{*j} & = & b_j, \qquad (j = 1, \dots, n).
					\label{columnsum2}
		\end{eqnarray}
Summing these two equations over $i$ and $j$ respectively and
combining with (\ref{zerosum}), we obtain
		\begin{equation}
\sum_j f_{*j} = 0.
		\end{equation}
 From property {\bf P6} it follows that
		\begin{equation}
f_{*j} = 0 \qquad (j = 1, \dots, n).
		\end{equation}
The forces are equal to zero on all edges adjacent to $P$.
Therefore $f_{ij}$ satisfies the constraints (\ref{rowsum}) to
(\ref{positivef}): it is a feasible solution.

 From property {\bf P5}, we derive
		\begin{equation}
f_{ij} \gamma_{ij} = 0.
		\end{equation}
Summing over $i$ and $j$, we obtain
		\begin{equation}
\sum_i a_i \alpha_i - \sum_j b_j \beta_j 
- \sum_i \sum_j f_{ij} c_{ij} = 0.
					\label{bidon2}
		\end{equation}
Consider another feasible solution $f'_{ij}$. From
(\ref{positivef}) and property {\bf P4} we have
		\begin{equation}
f'_{ij} \gamma_{ij} \ge 0
		\end{equation}
and therefore, summing over $i$ and $j$ and using (\ref{rowsum})
and (\ref{columnsum}):
		\begin{equation}
\sum_i a_i \alpha_i - \sum_j b_j \beta_j 
- \sum_i \sum_j f'_{ij} c_{ij} \ge 0.
		\end{equation}
Comparing with (\ref{bidon2}), we have
		\begin{equation}
\sum_i \sum_j f'_{ij} c_{ij} \le \sum_i \sum_j f_{ij} c_{ij}
		\end{equation}
which shows that $f_{ij}$ is an optimal solution.
\section{Notes on practical implementation}
\label{s:implementation}
\subsection{Dropping rows one by one}

Experience showed that the following modification of the algorithm
results in an important reduction in computing time (typically a
factor 3 for $m = n = 100$).  The stop $Q$ is cut into $m$
independent pieces $Q_1$, \dots, $Q_m$, each supporting one row,
and these stops are lowered one by one, each time waiting until an
equilibrium has been reached before starting the next stop. This
can be done by a simple modification of the algorithm described in
Section~\ref{s:formal}.
\subsection{General organization}

Measurements show that most of the computing time is spent in the
descent phase, and specifically in finding the minimum $d$ in
(\ref{d-def}). On the other hand, most of the complexity of the
program lies in updating the structure of the trees and the
associated information. Therefore only the computation of $d$
needs to be optimized for speed. In the remainder of the program,
one can freely use the structures which allow the easiest, most
natural, and most legible representation.

Experience shows that {\em triply linked trees}
(\cite[Section~2.3.3]{Knu73a}) are a convenient structure.  To
each row are associated three pointers to its mother, its eldest
daughter, and its next younger brother.  Similarly, to each column
are associated three pointers to its father, its eldest son, and
its next younger sister.  The first pointer is used to move
upwards in the tree, for instance in order to determine the
main path after contact has been made. The two other pointers are
used to explore a branch, for instance in order to change its
status from moving to fixed, or conversely, after a capture.

Each row has seven quantities associated with it: its weight $a_i$
(which does not change during the computation); its height
$\alpha_i$; the three pointers; the force of contact with its
mother; and a flag indicating whether the row belongs to the
moving tree or to the fixed tree.  Each column has seven similar
associated quantities.
\subsection{Methods for descent}
\label{s:descent}

Finding the minimum $d$ defined by (\ref{d-def}) would seem to be
a trivial task, involving two loops over $i$ and $j$ and about 15
lines of code. This would require a computing time of order $m n$
for each descent. In fact, there is room for considerable
improvement over this simple scheme.

\subsubsection{Version A}

First we observe that our task is to find the minimum among
quantities $\gamma_{ij}$ which are positive or zero. Therefore if,
during the examination of the $\gamma_{ij}$, we find one which is
zero, then we know that we have found the minimum and we can end
the search immediately. We will refer to the program incorporating
this simple device as {\em Version A}.

This is especially effective when the $c_{ij}$ take only a small
number of integer values. We will refer to this as the {\em
degenerate case} (see below, Section~\ref{s:hit-deg}). The
distances $\gamma_{ij}$ then take themselves only a small number
of distinct values. For large $m$ and $n$, each of these values
appears many times. In particular, as soon as the algorithm is in
progress, the value $\gamma_{ij} = 0$ typically appears many
times.  Therefore the search can be discontinued at an early time
and the computing time is much less than $O(m n)$.

When this method is used, experience shows that it is advisable to
start the search at a variable point in the $\gamma_{ij}$ matrix.
If the search is always started from the beginning, the contact
edge tends to be always selected from the same part of the matrix;
an unbalanced situation develops, and computing time increases.
\newpage
\vspace*{-10mm}
One method consists in choosing the starting point at random in
the matrix. This has the disadvantage of requiring the use of a
random number generator. A better solution (suggested by
A.~Noullez) is to compute the rank $r'$ of each new starting point
 from the rank $r$ of the previous one by a simple formula, such as
\begin{equation}
r' = r + \lfloor Kmn \rfloor \quad \pmod{m n},
\label{e:shuffling}
\end{equation}
where $[\cdot]$ denotes the integer part. 

The most uniform distribution of points is obtained by choosing $K
= (\sqrt5 - 1)/2$, the inverse of the golden ratio. This was found
to give very good results.

\subsubsection{Version B}

Another line of attack is based on the realization that the
structure of the system changes only partially from one cycle to
the next. Therefore we can try to save and re-use information on
the distances. (This approach was inspired by a study of
the LSAP algorithm presented in \cite[chap.~1]{BD80a} for the
particular case of the assignment problem.) In particular we may
try to take advantage of the tree structure which pervades the
algorithm, noting that much of that structure is left intact in a
capture episode. We will refer to the program developed along
these lines as {\em Version B}.

Several methods were tried. We describe here the method which gave
the best results , and which is incorporated in Version~B. We
define a {\em male branch} ${\cal A}_i$ as the branch (of the
moving or fixed tree) whose head is row $A_i$.  (Note that here is
a one-to-one correspondence between rows and male branches).  For
every pair $(i, j)$, we find the minimal distance $\delta_{ij}$
between the rows belonging to the male branch ${\cal A}_i$ and the
column $B_j$, and we note for which row this minimal distance is
realized. (Note that we consider here all male branches and all
columns, irrespective of whether they are moving or fixed).  When
the structure of the trees changes, this information is updated.
At the next descent, the updated information can then be used to
compute quickly the minimum $d$: if the moving stop is $Q_i$, the
whole moving tree consists of the male branch ${\cal A}_{i}$, and
one has only to find the minimum among the stored distances
$\delta_{ij}$ between that male branch and the fixed columns.
This takes a time $O(n)$.

The updating of the distances $\delta_{ij}$ is somewhat
complex. All male branches which have
their head on the main path need to be reconsidered. The cases of
capture by the fixed tree and by the moving tree require different
treatments. Also the three pieces of the main path determined by
the contact edge and the breaking edge have to be treated
separately. Savings in computing time are achieved by looking
for cases where it is not necessary to recompute the
$\delta_{ij}$.

Version B is more complex than Version A. It also requires about
twice as much memory, since the $\delta_{ij}$ array must be saved
in addition to the given $c_{ij}$ array. However, it is definitely
faster in the general, non-degenerate case (see below,
Section~\ref{s:hit}).
\section{Tests}
\label{s:tests}

Numerical tests were performed to verify the correctness of the
algorithm and to measure its performance. The computations were
made on a Hewlett-Packard Apollo Series 700, Model 720
workstation.

Comparisons were made with the subroutine \verb|H03ABF| of the NAG
library \cite{NAG91a}.  In all computed cases (which number in the
thousands) it was verified that exactly the same optimal cost
(\ref{maxisum}) is found with the present algorithm and with the
NAG algorithm.

The values $a_i$ and $b_j$ are taken as positive integers. They
are first chosen at random in the intervals
\begin{equation}
1 \le a_i \le a_{\rm max}, \qquad 1 \le b_j \le b_{\rm max},
\end{equation}
where $a_{\rm max}$ and $b_{\rm max}$ are two constants satisfying
$m a_{\rm max} = n b_{\rm max}$. Small adjustments are then made
in order to satisfy the relation (\ref{zerosum}) exactly.

Tests show that the computing time is insensitive
to the values of $a_{\rm max}$ and $b_{\rm max}$, provided that
they are not too close to unity (A variation becomes
detectable for values of 10 or less). In practice we
take $a_{\rm max} = 160000 / m$, $b_{\rm max} = 160000 / n$.

The values $c_{ij}$ are also taken as integers, randomly chosen in
the interval
\begin{equation}
1 \le c_{ij} \le c_{\rm max}.
					\label{c-range}
\end{equation}
Again the computing time is found to be insensitive to the value
of $c_{\rm max}$, provided that it is large enough. Tests show
that the relevant quantity is the ratio
\begin{equation}
c_* = {c_{\rm max} \over m n}.
\end{equation}
Variations of the computing time begin to be noticeable when $c_*$
is less than 1 (see below Section~\ref{s:hit-deg}). This
corresponds to the onset of degeneracy: for $c_* \ll 1$, each
value in the allowed range (\ref{c-range}) appears many times in
the $c_{ij}$ matrix. Thus, the value taken for $c_{\rm max}$
depends on whether the non-degenerate or the degenerate case is
considered.

For simplicity only the square case $m = n$ was considered, with
$n$ ranging from 10 to 1000. Time measurements were generally
averaged over a series of 100 computations, in order to obtain
more accurate values.
\subsection{Non-degenerate case}
\label{s:hit}

We take $c_* = 1$. Note that this corresponds to a case where each
value in the allowed interval (\ref{c-range}) is present
once on average in the matrix. 

\begin{figure}[hbtp]
\centerline{\psfig{file=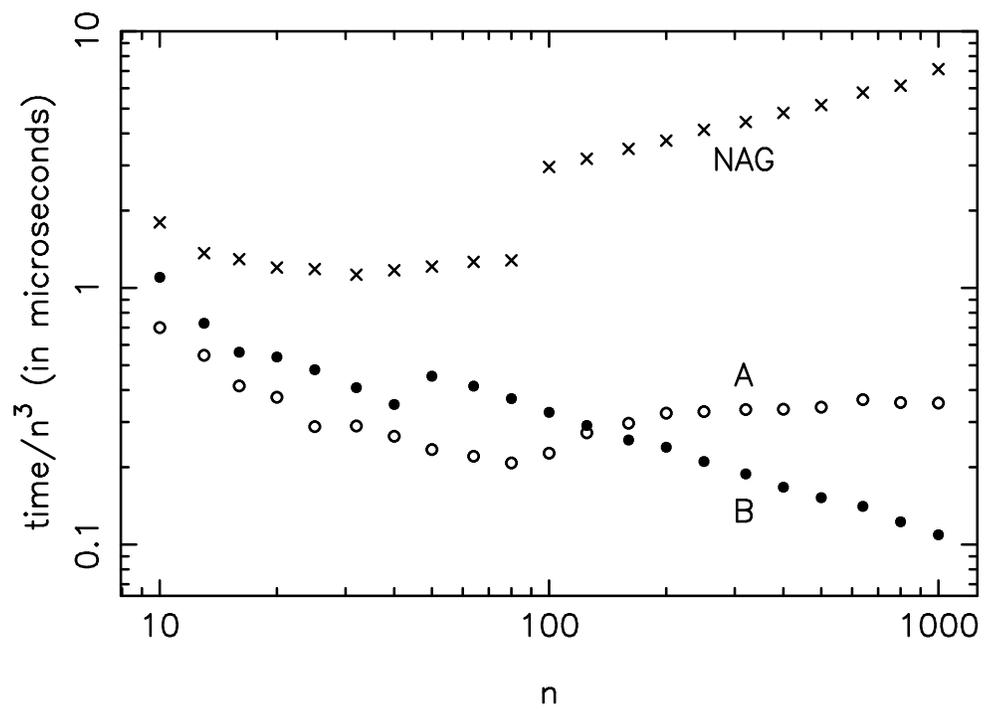,width=13cm,%
      bbllx=28mm,bblly=40mm,bburx=173mm,bbury=143mm}}
\caption{transportation problem: computing time (divided by $n^3$) as a function
of size $n$. Crosses: NAG subroutine. Open circles: present
algorithm, Version A. Filled circles: present algorithm, Version B.}
	\label{f:hit}	\end{figure}

Fig.~\ref{f:hit} shows the computing time (divided by $n^3$ for
better clarity) as a function of $n$, for three algorithms:

\begin{itemize}

\item Crosses correspond to the subroutine \verb|H03ABF| of the NAG
library \cite{NAG91a}. The time appears to grow asymptotically as
$n^{3.35}$. (A curious discontinuity is present: the computing
time jumps up suddenly by a factor of about 2.2 between the
values $n = 94$ and $n = 95$. This is probably due to
peculiarities of the computer hardware.)

\item Open circles represent Version A of the present algorithm
(see Section~\ref{s:descent}). Computing time grows asymptotically
as $n^{3.05}$.

\item Filled circles represent Version B of the present algorithm
(see Section~\ref{s:descent}). For low values of $n$, the
computation is slower than with Version A because of the extra work
involved in computing the distances $\delta_{ij}$. Above $n =
100$, however, this extra work begins to pay
off. Computing time grows asymptotically as $n^{2.5}$.

\end{itemize}

Version B is clearly the best method. For a $1000 \times 1000$
problem, the NAG subroutine takes about 7000 seconds, while
Version B takes about 110 seconds. The ratio increases for larger
values of $n$.
\subsection{Degenerate case}
\label{s:hit-deg}

A value $c_{\rm max} = 20$ was chosen as representative for a
degenerate problem.  In particular, this is a typical value for
applications to lattice gas problems \cite{Hen89a}. Thus, $c_{ij}$
can take only integer values from 1 to 20.

\begin{figure}[hbtp]
\centerline{\psfig{file=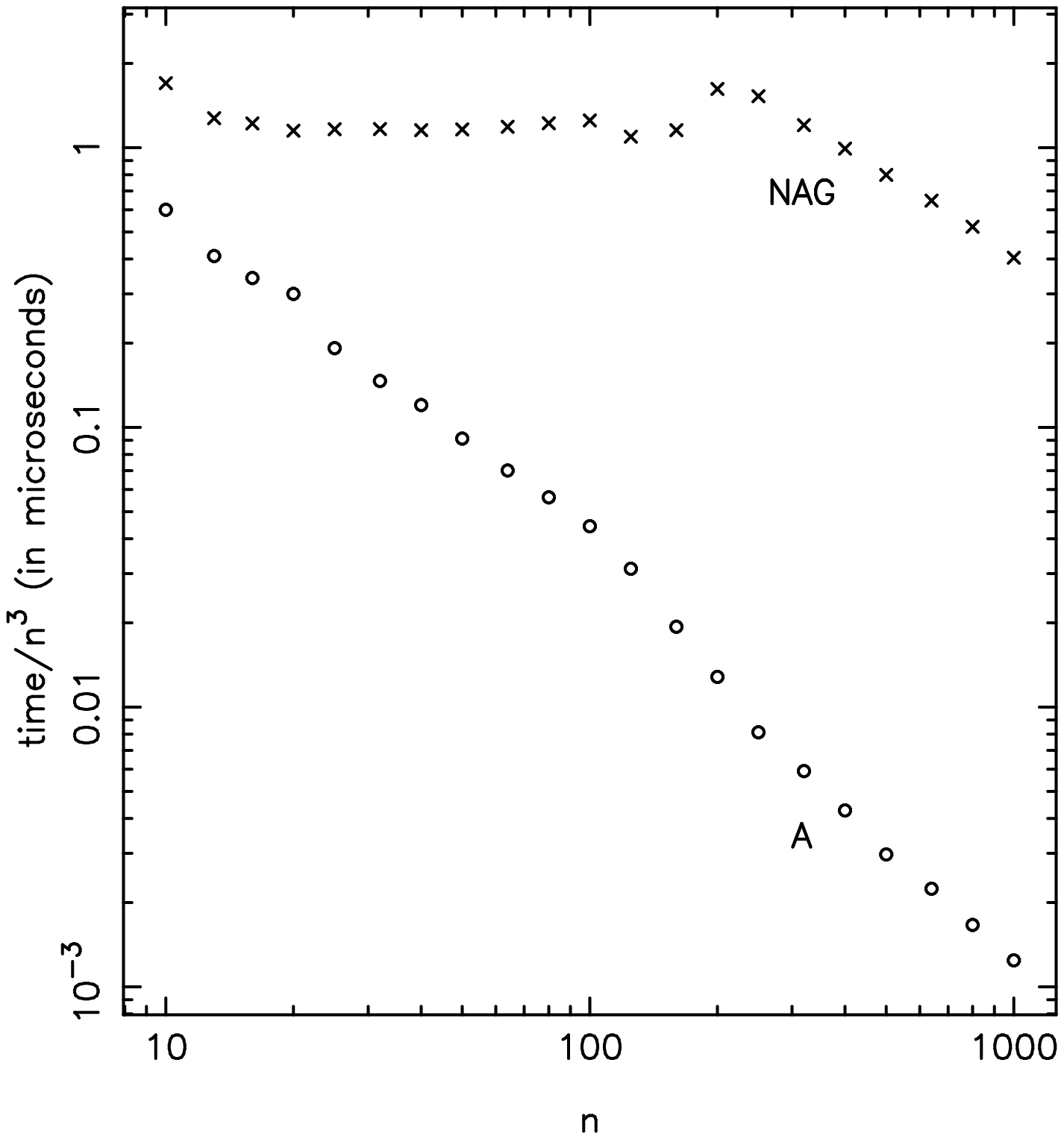,width=13cm,%
      bbllx=28mm,bblly=15mm,bburx=172mm,bbury=169mm}}
\caption{Degenerate transportation problem: computing time (divided by $n^3$) as
a function of size $n$. Crosses: NAG subroutine. Open circles:
present algorithm, Version A.}
	\label{f:hit_deg}
	\end{figure}

Fig.~\ref{f:hit_deg} shows computing times as a function of $n$,
for two algorithms: the NAG subroutine and Version A of the present
algorithm. (Version~B is inefficient in the degenerate case and is
not shown). 

\begin{itemize}

\item For the NAG subroutine (crosses), computing time is
essentially the same as in the non-degenerate case
(Fig.~\ref{f:hit}) up to about $n = 100$. For larger values, the
effect of the degeneracy begins to be felt and the slope
decreases.  Asymptotically, the computing time appears to grow as
about $n^{2}$.

\item For the present algorithm (open circles), the decrease in
computing time with respect to the non-degenerate case is more
marked and starts earlier, at about $n = 20$. The time dependence
is more complex.  The final slope indicates a dependence in
$n^{1.65}$.

\end{itemize}

For a $1000 \times 1000$ problem, the NAG subroutine takes 400
seconds, while Version A takes about 1.25 seconds. The ratio again
increases for larger values of $n$.

We remark that an exponent of $n$ less than 2 means that for large
values of $n$, the time needed to solve the problem is small
compared to the time needed to set it up, since simply copying the
$c_{ij}$ matrix into memory takes a time proportional to $n^2$ !
Note also that in this situation, most $c_{ij}$ values will never
be used.

\clearpage
\subsection{Assignment problem}
\label{s:ass}

Tests were also made for the particular case of the assignment
problem, where $m = n$ and all rods have weights $a_i = 1$, $b_j =
1$. (It is then easily shown that an optimal solution obtained
with the present algorithm automatically satisfies the constraint
(\ref{0or1}). In that case, comparisons can also be made with the
subroutine {\tt LSAP} of Burkard and Derigs \cite{BD80a} (noted BD
here).

\begin{figure}[hbtp]
\centerline{\psfig{file=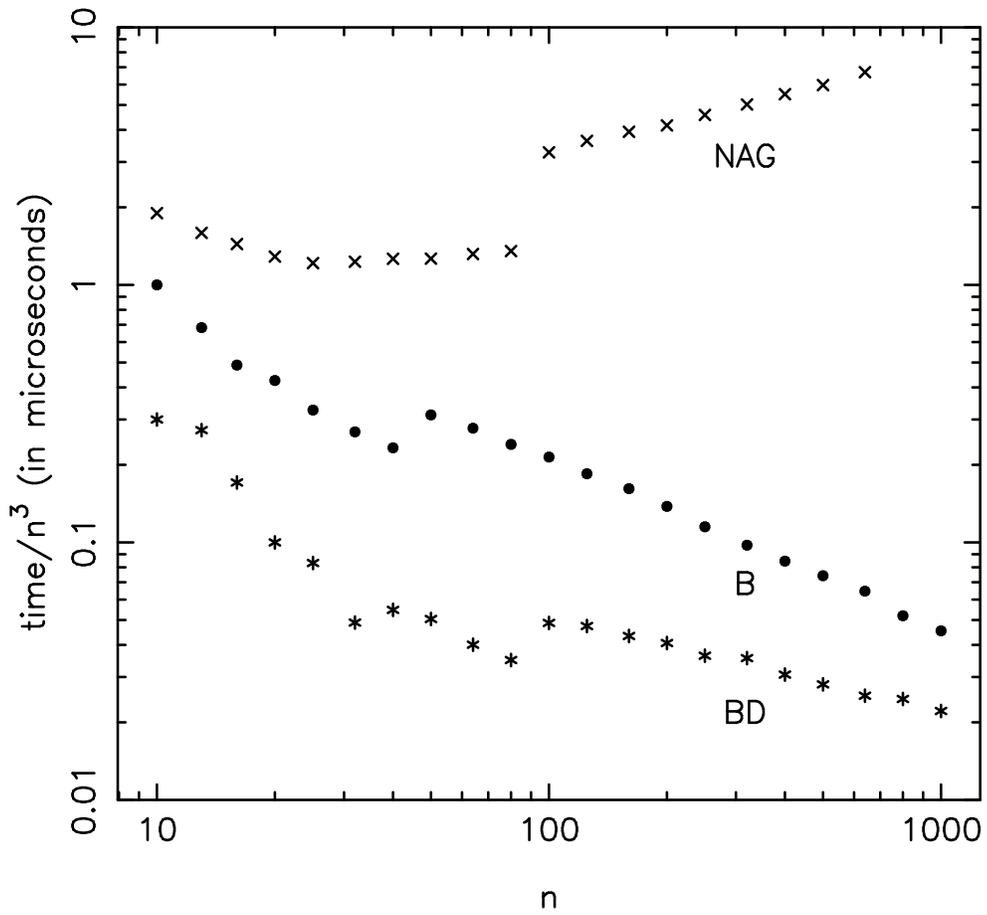,width=13cm,%
      bbllx=28mm,bblly=66mm,bburx=172mm,bbury=199mm}}
\caption{Assignment problem: computing time (divided by $n^3$) as a
function of size $n$. Crosses: NAG subroutine. Filled circles:
present algorithm, Version~B. Asterisks: Burkard and Derigs
algorithm.}
	\label{f:ass}	\end{figure}

Fig.~\ref{f:ass} compares computing times in the general
(non-degenerate) case with $c_* = 1$. As before, the NAG
subroutine is much slower. The BD algorithm is fastest: for a
$1000 \times 1000$ problem, computing time is about 45 seconds
for Version~B and 22 seconds for the BD algorithm. The difference
decreases when $n$ increases, however: the asymptotic law is about
$n^{2.4}$ for Version~B, compared to $n^{2.7}$ for BD.

\begin{figure}[hbtp]
\centerline{\psfig{file=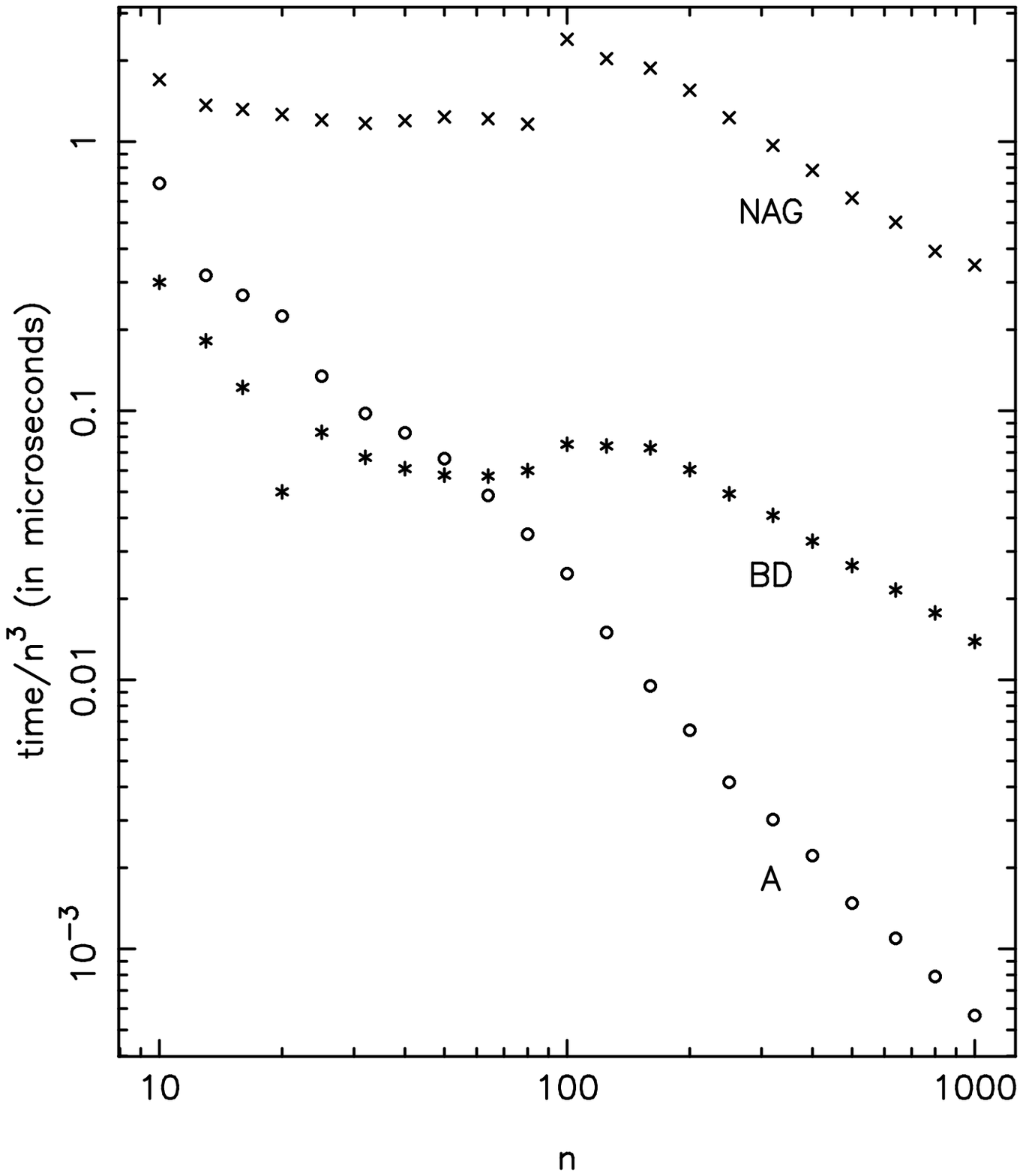,width=13cm,%
      bbllx=28mm,bblly=15mm,bburx=173mm,bbury=181mm}}
\caption{Degenerate assignment problem: computing time (divided by $n^3$)
as a function of size $n$. Crosses: NAG subroutine. Open circles:
present algorithm, Version~A. Asterisks: Burkard and Derigs
algorithm.}
	\label{f:ass_deg}
	\end{figure}

Fig.~\ref{f:ass_deg} compares NAG, Version~A, and BD for the
degenerate assignment problem, with $c_{\rm max} = 20$. Here
Version~A is fastest: for a $1000
\times 1000$ problem, computing time is 0.6 seconds for Version~A,
and 14 seconds for BD. The difference increases with $n$: the
asymptotic behaviour is in $n^{1.55}$ for Version~A, $n^2$ for BD
and NAG.

\section{Final comments}
\label{s:conclusions}

1. Many variations of the present algorithm could be imagined, and
it is quite possible that some of them would increase its speed.
For instance, sometimes an arbitrary choice can be made for the
contact edge or the breaking edge (see steps 3 and 5 in
Section~\ref{s:formal}); one might try to determine what is the
best choice.
 
2. As mentioned in the Introduction, the present algorithm was
conceived and developed independently, without recourse to the
existing litterature. Looking back, however, it becomes clear that
some relations exist.  The heights of the rods, for instance, are
nothing else than the classical {\em dual variables}; hence the
choice of the customary notations $\alpha_i$ and $\beta_j$ in
the present paper.

In particular, after this work was completed, M.~Hartmann called
our attention to the references \cite{SJB82a} and \cite{ABG84a}.
These papers describe algorithms for the minimum cost flow
problem, which includes the transportation problem as a special case.
One starts from zero flow and continually increments it,
maintaining at all times a minimum-cost solution, until the
desired flow is attained. It appears that the algorithm of the
present paper belongs essentially to the same family, known as
{\em parametric algorithms}. The equivalent of the flow in the
present model is the total force applied by the lines to the
columns; as is easily shown, this total force starts from zero and
increases with time, until it is equal to the total weight of the
lines. More specifically, it increases by $\lambda$ during each
readjustment of the forces (Section~\ref{s:algorithm}, Step 4).

3. As mentioned in Section~\ref{s:implementation}, most of the
computing time in the present algorithm is spent on finding the
minimum in a large set of numbers. Efficient algorithms have been
developed for this operation on massively parallel computers
\cite{HS86a}. A parallel implementation of the algorithm might
therefore be of interest.

4. A comparison of figures \ref{f:hit} and \ref{f:hit_deg}, or
\ref{f:ass} and \ref{f:ass_deg}, shows that the solution of the
degenerate problem is much faster.  This might be of interest in
situations where a great accuracy is not required, or is not
present in the data. In such cases, it will be very advantageous
to round off the $c_{ij}$ values so as to reduce them to a
comparatively small set of values.

5. The idea of using mechanical analog computers for optimization
problems is not new. For instance, reference \cite{VSD86a}
describes a mechanical device made of shafts and gears, which
can in principle solve the general instance of the linear
programming problem. However, this device is introduced in
\cite{VSD86a} only as a conceptual tool in a theoretical study of 
the complexity of analog computation; it is not intended as a
model for a practical algorithm. We remark also that the
mechanical model of the present paper is adapted to the special
case of the transportation problem, and is therefore much simpler (and
presumably more efficient) in that special case. As a rough
measure, the present model has $O(m+n)$ moving parts, while the
model of \cite{VSD86a} would have $O(mn)$ moving parts in a $m
\times n$ transportation problem.

In this connection, it is natural to ask whether the mechanical
model used here to simulate the transportation problem can be
extended to the more general minimum cost flow problem, or to the
even more general linear programming problem. We have not found
any obvious way to do this.
\appendix
\section{Bounds on the maximal number of cycles}
\label{s:number-of-cycles}

We call $Z(m, n)$ the maximal number of cycles, for the $m
\times n$ problem, using the algorithm described in
Section~\ref{s:formal}.  We derive here some rigorous bounds on
this number.
\subsection{Upper bound}

We derive first a general upper bound for $Z(m, n)$.

We number with an index $l$ the successive levels of the moving
tree. The stop $Q$ is at level $l = 0$, the sons of $Q$ are at
level $l = 1$, the grand-daughters of $Q$ are at level $l = 2$,
and so on. Note that odd levels correspond to rows and even
levels to columns. We call $g_l$ the number of nodes of the tree
at level $l$. The sequence of numbers $g_1$, $g_2$, \dots, will be
called the {\em signature} of the moving tree.

We consider now all possible signatures for given $m$ and $n$,
assuming that the final state has not yet been reached, i.e. that
the moving tree still contains at least one row and the fixed
tree still contains at least one column. We order these signatures
as follows. First we sort by decreasing $g_1$. Next we sort
each subset by increasing $g_2$. Next we sort each subsubset
(corresponding to given $g_1$ and $g_2$) by decreasing $g_3$; and
so on, always using decreasing order for odd values of $l$ and
increasing order for even values. Finally, we number the sorted
signatures with $K = 1$, 2, \dots.

It is then easy to show that {\em $K$ always increases during a
cycle}. There are two cases:

\begin{enumerate}

\item The moving tree captures a subtree from the fixed tree. The
root of this subtree (after the capture) is a column.
Therefore the first level $l$ at which there is a change in the
signature is even, and $g_l$ increases by one unit. From the above
sorting method it follows that $K$ increases.

\item The moving tree loses a subtree. The root of this subtree
(before the capture) is a row. Therefore the first level $l$ at
which there is a change in the signature is odd, and $g_l$
decreases by one unit. Again $K$ increases.

\end{enumerate}

Therefore we obtain an upper bound on the number of cycles simply
by counting the signatures. We call $p = g_1 + g_3 + \dots$ the
number of rows, and $q = g_2 + g_4 + \dots$ the number of columns
in the moving tree. Since the moving tree is assumed to be
non-empty, $p$ can take values from 1 to $m$. Similarly, since the
fixed tree is non-empty, $q$ can take values from 0 to $n - 1$. We
evaluate first the number of signatures for given $p$ and $q$. A
signature can also be represented by a sequence of $p + q$ binary
digits: we write $g_1$ digits 1, then $g_2$ digits 0, then $g_3$
digits 1, and so on. The first digit must be a 1. There are $q$
digits 0, which can be placed anywhere in the remaining $p + q -
1$ positions.  Therefore the number of possible signatures is
\begin{equation}
{p + q - 1 \choose q}.
\end{equation}
Summing over $p$ and $q$, we obtain the following upper bound for
$Z$:
\begin{equation}
Z(m,n) \le Z_{\rm sup}(m,n) = {m + n \choose m} - 1.
					\label{upper-bound}
\end{equation}

Exactly the same considerations can be applied to the fixed tree.
 We number with $l$ the successive levels. We call $g'_l$ the
 number of nodes at level $l$. The sequence of numbers $g'_1$,
 $g'_2$, \dots, will be called the signature of the fixed
tree. We sort the signatures, using decreasing order for odd $l$
(columns) and increasing order for even $l$ (rows). We number the
sorted signatures with $K' = 1, 2, \dots$. The number $K'$ also
always increases during a cycle. We call $q' = g'_1 + g'_3 +
\dots$ the number of columns, and $p' = g'_2 + g'_4 + \dots$ the
number of rows in the fixed tree. Counting the number of
signatures as above, we find

\begin{equation}
{m + n \choose n} - 1
\end{equation}
i.e. exactly the same upper bound as in (\ref{upper-bound}).

\subsection{Better upper bound}

A much better upper bound can be obtained by noting that only some
combinations of $K$ and $K'$ are permitted. This is because we
must have
\begin{equation}
p + p' = m, \qquad q + q' = n.
\end{equation}
Thus, in a $(K, K')$ plane, only a subset of points are allowed.
Combining this with the fact that both $K$ and $K'$ must increase
at each step, one can trace the possible paths in the plane and
derive an upper limit on the number of steps, which we call
$Z'_{\rm sup}(m,n)$.

Unfortunately a general formula giving $Z'_{\rm sup}(m,n)$ for
arbitrary $m$ and $n$ has not been found. Results obtained by a
computer program for values of $m$ and $n$ up to 10 are listed in
Table~\ref{t:K-and-K'}.

		\begin{table}[hbt]

\caption{Upper limit $Z'_{\rm sup}(m,n)$ on the number of cycles
for $1 \le m \le 10$, $1 \le n \le 10$.}

\label{t:K-and-K'}
		$$\begin{tabular}{|c|r|rrrrrrrrrr|}
\hline
\multicolumn{2}{|c|}{} & \multicolumn{10}{c|}{$n$} \\
\cline{3-12}
\multicolumn{2}{|c|}{} & 1 & 2 & 3 & 4 & 5 & 6 & 7 & 8 & 9 & 10 \\
\hline
& 1 &   1 &   2 &   3 &   4 &   5 &   6 &   7 &   8 &   9 &  10 \\
& 2 &   2 &   4 &   6 &   8 &  10 &  12 &  14 &  16 &  18 &  20 \\
& 3 &   3 &   6 &  10 &  14 &  19 &  24 &  30 &  36 &  43 &  50 \\
& 4 &   4 &   8 &  14 &  22 &  30 &  40 &  52 &  64 &  78 &  94 \\
$m$&5&  5 &  10 &  19 &  30 &  46 &  62 &  83 & 108 & 138 & 170 \\
& 6 &   6 &  12 &  24 &  40 &  62 &  94 & 126 & 168 & 222 & 284 \\
& 7 &   7 &  14 &  30 &  52 &  83 & 126 & 190 & 254 & 339 & 448 \\
& 8 &   8 &  16 &  36 &  64 & 108 & 168 & 254 & 382 & 510 & 682 \\
& 9 &   9 &  18 &  43 &  78 & 138 & 222 & 339 & 510 & 766 &1022 \\
&10 &  10 &  20 &  50 &  94 & 170 & 284 & 448 & 682 &1022 &1534 \\
\hline
		\end{tabular}$$
 		\end{table} 

For $m = n = 10$, for instance, we have $Z_{\rm sup}(m,n) =
184755$,  $Z'_{\rm sup}(m,n) = 1534$.
\subsection{Lower bound}

In the square case $m = n$, the following lower bound can be
proved:

\begin{equation}
Z(n,n) \ge Z_{\rm inf}(n,n) = 3 \times 2^{n - 1} - 2.
					\label{lower-bound}
\end{equation}

The proof is cumbersome and will not be given here. It consists in
showing that the number of cycles equals $Z_{\rm inf}$ in the
following case:

\begin{eqnarray}
a_1 & = & 1, \qquad b_1 = 2, \nonumber \\
a_i & = & a_{i-1} + b_{i-1}, \qquad b_i = b_{i-1} + a_i, 
  \qquad (i = 2, \dots, n - 1), \nonumber \\
a_n & = & a_{n-1} + b_{n-1}, \qquad b_n = b_{n-1}, \nonumber \\
c_{ij} & = & \cases{
(n + 1 - i)(n + 1 - j) - n^2(2^i + 2^j) & if $i = j$, \cr
(n + 1 - i)(n + 1 - j) - n^2(2^j) & if $i < j$, \cr
(n + 1 - i)(n + 1 - j) - n^2(2^i) & if $i > j$. \cr
}
					\label{lower-example}
\end{eqnarray}
This was also verified by a direct
application of the numerical algorithm for $n = 1$ to 17.
Note that the sequence $a_1$, $b_1$, $a_2$, $b_2$, \dots is 
the Fibonacci sequence, minus its first term.

Comparing with the diagonal of Table~\ref{t:K-and-K'}, we find
that the upper bound given in that Table is identical to the lower
bound given by (\ref{lower-bound}). Therefore, for $n = 1$ to 10,
we know the exact value of the maximal number of cycles, which is

\begin{equation}
Z(n,n) = 3 \times 2^{n - 1} - 2,
					\label{exact-Z}
\end{equation}
and (\ref{lower-example}) is a worst case, achieving this maximal
value.  There is a strong suggestion that (\ref{exact-Z}) holds
for all values of $n$, but this has not been proved.
\subsection{Comparison with observed values}

We observe that (\ref{exact-Z}) corresponds to a computing time
which grows exponentially with $n$. Fortunately, numerical tests
with randomly chosen examples show a much milder increase, which
is approximately linear in~$n$.  Table~\ref{t:comp} compares the
values (\ref{exact-Z}) with the average observed values of the
number of cycles, for Version B of the algorithm, in the
non-degenerate case.  The r.m.s. dispersions are also given; they
show that individual values do not deviate much from the average.

		\begin{table}[hbt]
\caption{Number of cycles in the ``bad'' case
(\protect\ref{lower-example}) (column 2), and observed
number of cycles (column 3).}
\label{t:comp}
		$$\begin{tabular}{|r|r|r|}
\hline
$n$ & $Z_{\rm inf}(n,n)$ & observed \\
\hline
  5 &         46 & $ 10 \pm  1$ \\
 10 &       1534 & $ 23 \pm  2$\\
 15 &      49150 & $ 37 \pm  3$\\
 20 &    1572862 & $ 50 \pm  4$\\
 25 &   50331646 & $ 64 \pm  5$\\
 30 & 1610612734 & $ 80 \pm  6$\\
 50 &            & $141 \pm 10$\\
100 &            & $304 \pm 17$\\
\hline
		\end{tabular}$$
		\end{table}

This considerable difference between the ``bad case''
(\ref{lower-example})and the average case can be probably
understood by noting that (\ref{lower-example}) is a rather
extreme case: the coefficients $a_i$, $b_j$, $c_{ij}$ form
essentially geometrical progressions. For $m = 1000$, for
instance, the ratio $a_{1000}/a_1$ is of the order of $10^{400}$.
This is not likely to be encountered in applications.
\subsection{Acknowledgements}

I thank P.~Bernhard, U.~Frisch, M.~Hartmann, J.~Morgenstern,
A.~Noullez, K.~Steiglitz, and S.~Stidham for discussions and
comments.

\newpage
\section{Addendum (September 2002)}

The present paper was submitted to {\em Mathematical Programming}
in May 1992. The following referee's report was subsequently
received:

\vspace{5mm}
\noindent {\em (Beginning of referee's report)}
\vspace{0mm}

\subsubsection*{General comments}

This well-written paper proposes an algorithm for the
transportation problem that is motivated by an analog model, and
has some comparisons of a code for this algorithm with a NAG
code. The new algorithm is found to be much faster than the NAG
code on randomly generated problems.

The author of this paper comes from outside the Math Programming
community, and is to be commended for taking the time to bring
his/her fresh perspective to the subject. However, it is tempting
to dismiss this paper on the basis that the author has merely
re-discovered a known algorithm, namely the one in [13], noting
that [13] already contains a computational comparison of such an
algorithm with other algorithms existing at the time. A further
reasonable criticism of the current paper is that the NAG routine
``H03ABF uses the `stepping stone' method, modified to accept
degenerate cases'' (quoting from the NAG documentation),
referenced in the book ``Linear Programming'' by G. Hadley, 1962
(this description should have been in the paper). The routine
itself has been in the NAG library at least from 1975. Thus this
code does not represent the current state of the art in
transportation algorithms, as would be found, e.g., in Ahuja et
al. Also, comparing codes on randomly generated problems can be
misleading in any case.

On the other hand, I believe that in general such outside
contributions should be valued, and specifically that this paper
has something to offer if it is drastically re-written:

\begin{enumerate}

\item The application of transportation problems to lattice gas
models is intriguing, and is worth further discussion. Is there a
short way to say how such models arise and the significance of the
transportation subproblem? Also, do such problems tend to be
sparse or dense, how big are $m$ and $n$ in typical problems, how
large do the supplies and demands tend to be, and how large do the
costs tend to be? These are all parameters that are important in
assessing which of the standard modern algorithms might work well
on such problems.

\item The computational testing in the paper needs to be fixed.
It would be more believable if the codes were tested on instances
arising from the actual application rather than random
problems. It would also be a service to the community to release
several typical examples of such problems to, e.g., the DIMACS
library of network flow instances, so that other codes can be
tried on these problems. DIMACS also has available some
well-tested random network generators; it enhances comparability
between codes if standard generators are used instead of {\em ad
hoc} generators. DIMACS codes and generators are available via
anonymous ftp from dimacs.rutgers.edu in directory pub/netflow.

\item I enjoyed the analog model and the physical insight it gives
to the algorithm. In fact, I believe that a slightly different
presentation of the algorithm that ties it more closely to the
physical model would improve the paper: My suggestion is to make
$P$ and $Q$ into an extra row and column (which they effectively
are in the model anyway), both with zero weight/buoyancy, and such
that $c_{Pj} = c_{iQ} = 0$ for all $i$ and $j$. However, to act as
stops in the initial configuration, we must give $P$ an artificial
weight of $\sum b_j$, and $Q$ an initial buoyancy of $\sum
a_i$. The initial configuration is in fact optimal with these
artificial weights. The aim of the algorithm is then to decrease
the artificial weights to zero while maintaining optimality at
each intermediate artificial weight value. Thinking of $P$ and $Q$
as lines remove the need to treat them as special cases elsewhere.

\item Much of the proof of the algorithm's correctness is
unnecessary since the author is just rediscovering well-known
arguments, and has proposed an algorithm that fits nicely into
known classes of algorithms. For example, equation (25) is known
as complementary slackness, and the fact that complementary
slackness plus feasibility equals optimality is so well-known that
it can be stated without a reference. Although the author derived
the algorithm independently of known algorithms, in order to
effectively present it to an audience which is familiar with known
algorithms it would be helpful to discuss the algorithm as if it
were a special case of what is called ``dual node-infeasible''
simplex algorithms in [13]. I have in mind something like the
following (assuming that suggestion 3 above is taken): The
algorithm maintains dual feasible variables (the heights) and
primal flows (forces; this correspondence between heights/forces
and dual/primal variables is too important to leave it to a
comment in the conclusion) that satisfy non-negativity and all
supplies and demands except possibly at the extra nodes $P$ and
$Q$, and ensures that the primal and dual variables satisfy
complementary slackness. It also maintains a basic tree, namely
the fixed tree plus the moving tree, plus the extra arc where
fresh contact between the fixed and moving trees occurs. The extra
arc allows some of the surplus supply at $P$ to be pushed through
the tree to cancel out some of the surplus demand at $Q$. An arc
whose flow drops to zero during this flow push can then be the
dual simplex outarc, and a standard dual simplex pivot (whose two
sides will be the fixed tree and the moving tree) will determine
which is the new extra inarc to be added to the tree.  This
continues until there is no surplus supply or demand at $P$ or
$Q$.  This change would allow Section 6 to be reduced to a few
sentences in Section 4, since it is well-understood that
node-infeasible dual network simplex maintains dual feasibility,
primal feasibility except for conservation, and complementary
slackness, and that when conservation is achieved we must be
optimal. Without this change and its attendant severe shortening,
the paper is not publishable since large parts of it recapitulate
familiar arguments.

\item The Appendix should be shortened and moved into the text to
establish what is known about how fast the algorithm
converges. All that is needed is the observation that the
algorithm is always finite without the need for any
anti-degeneracy device (since signatures are lexicographically
strictly increasing), that the upper bound $\scriptsize \left(
\begin{array} {c} m + n \\ n \end{array} \right)$ is easy to
derive, and that $n \times n$ examples exist which use $O(2^n)$
pivots. It should be pointed out that each piece of data in these
$n \times n$ examples has only $O(n)$ bits, so the examples show
that the algorithm is not even weakly polynomial. Also, either
``iterations'' or ``pivots'' is preferable to the word "cycles",
since ``cycles'' is too suggestive of cycles in graphs.

\item It seems possible that by using scaling (see, e.g., Ahuja et
al.), a weakly polynomial version of the algorithm could be
developed. The rough idea would be to scale the supplies and
demands.  Although in general this would lead to a problem where
$\sum a_i \ne \sum b_j$, the extra nodes would handily absorb the
extra flow. When we optimize at one scale factor and want to move
to the next one, supply or demand at a node might increase by one,
and we need to regain optimality.  I believe that we can do this
by a shortest path computation at each such node, looking for the
shortest distance path that will allow us to rehang the extra flow
from the proper one of $P$ or $Q$. This can be done in polynomial
time. We then have a problem where the total surplus supply/demand
at $P$ and $Q$ is $O(m+n)$. Degeneracy could cause a problem here
since it makes it difficult to get a polynomial bound on
iterations before optimality. I believe that the same shortest
path trick can be used to ensure that each iteration moves at
least one unit of flow. This shortest path business compresses
several iterations into a single iteration, and might be useful in
general. (Ideally, it would be nice to see a computational
comparison of the shortest path version with the tree version of
the algorithm on degenerate problems.) In addition, the shortest
path neatly ties in with the notion of the reduced costs as
``distances'', and is defensible in terms of the model as looking
for the smallest distance to move the moving tree so that some of
the force on $P$ and $Q$ can be lessened. This form of the
algorithm starts to look very much like the well-known successive
shortest paths algorithm (see Ahuja et al.).  In any case, a
mention of a possible scaling version of the algorithm would be
useful.

\item Most readers will prefer the term ``transportation problem''
to ``Hitchcock problem''.

\end{enumerate}

\vspace{1mm}
\noindent {\em (End of referee's report)}
\vspace{5mm}

Incidentally, I wish to thank here the unknown referee for
taking the trouble to write such a detailed and helpful report.

It took me some time to obtain the necessary papers, to assimilate
them and to re-do the calculations. In the end, I was able to
convince myself that, unfortunately, the referee was entirely
correct:

(i) The algorithm is not new. It is in fact identical, apart from
some trivial changes, with the ``parametric network simplex
algorithm'' described for instance in the book ``Network flows'',
published in 1993 by Ahuja, Magnanti and Orlin, pages 433--437.

(ii) The NAG algorithm which I used as a basis for comparison is
not state of the art. I made comparisons with more modern
algorithms, specifically with the results of Jianxiu Hao and
George Kocur in the paper ``An implementation of a shortest
augmenting path algorithm for the assignment problem'' (found on
the Internet server {\tt dimacs.rutgers.edu}, dated 1992), and my
program does not show a marked advantage anymore.

Unfortunately I did not feel able to follow the suggestion of the
referee and to drastically rewrite the paper; it would have been
too much work, in a field with which I am not very familiar. So I
published a much shorter version, including only the mechanical
model (Section~3 of the present paper) in {\em Comptes Rendus de
l'Acad\'emie des Sciences, Paris}, {\bf 321}, S\'erie I, 741--745
(1995).

Recently, some of my colleagues have used ideas from the present
paper to solve a problem in cosmology (see ``A reconstruction of
the initial conditions of the Universe by optimal mass
transportation'', by Uriel Frisch, Sabino Matarrese, Roya
Mohayaee, Andrei Sobolevski, {\em Nature}, {\bf 417}, 260--262
(2002) = astro-ph/0109483), and they asked me to make the
paper generally available by submitting it to arXiv.

The present text is unchanged from the 1992 original, with two
exceptions: a minor error has been corrected in Equ.~(47), and the
expression {\em Hitchcock problem} has been replaced by the more
modern name {\em transportation problem}, as suggested by the
referee.

\end{document}